\def\LaTeX{%
  \let\Begin\begin
  \let\End\end
  \let\salta\relax
  \let\finqui\relax
  \let\futuro\relax}

\def\UK{\def\our{our}\let\sz s}
\def\USA{\def\our{or}\let\sz z}

\UK
\LaTeX

\USA

\documentclass[twoside,12pt]{article}
\salta

\usepackage[usenames,dvipsnames]{color}
\usepackage{amsmath}
\usepackage{amsthm}
\usepackage{amssymb}
\usepackage[mathcal]{euscript}
\usepackage{cite}
\definecolor{viola}{rgb}{0.3,0,0.7}
\definecolor{ciclamino}{rgb}{0.5,0,0.5}
\definecolor{rosso}{rgb}{0.8,0,0}

\def\gianni #1{#1}
\def\pier #1{#1}
\def\betti #1{#1}
\def\gabri#1{#1}
\def\pcol #1{#1}

\finqui

\def\Beq{\Begin{equation}}
\def\Eeq{\End{equation}}
\def\Bsist{\Begin{eqnarray}}
\def\Esist{\End{eqnarray}}

\def\Bthm{\Begin{theorem}}
\def\Ethm{\End{theorem}}

\def\Bprop{\Begin{proposition}}
\def\Eprop{\End{proposition}}
\def\Bcor{\Begin{corollary}}
\def\Ecor{\End{corollary}}
\def\Brem{\Begin{remark}\rm}
\def\Erem{\End{remark}}

\def\Bdim{\Begin{proof}}
\def\Edim{\End{proof}}

\let\non\nonumber

\def\step #1 \par{\medskip\noindent{\bf #1.}\quad}

\def\Lip{Lip\-schitz}
\def\Holder{H\"older}
\def\Frechet{Fr\'echet}
\def\aand{\quad\hbox{and}\quad}

\def\wk{well-known}

\def\lhs{left-hand side}
\def\rhs{right-hand side}

\def\generaliz{generali\sz}
\def\lineariz{lineari\sz}

\def\regulariz{regulari\sz}

\def\bhv{behavi\our}

\def\multibold #1{\def\arg{#1}%
  \ifx\arg\pto \let\next\relax
  \else
  \def\next{\expandafter
    \def\csname #1#1#1\endcsname{{\bf #1}}%
    \multibold}%
  \fi \next}

\def\pto{.}

\def\multical #1{\def\arg{#1}%
  \ifx\arg\pto \let\next\relax
  \else
  \def\next{\expandafter
    \def\csname cal#1\endcsname{{\cal #1}}%
    \multical}%
  \fi \next}

\def\multimathop #1 {\def\arg{#1}%
  \ifx\arg\pto \let\next\relax
  \else
  \def\next{\expandafter
    \def\csname #1\endcsname{\mathop{\rm #1}\nolimits}%
    \multimathop}%
  \fi \next}

\multibold
qwertyuiopasdfghjklzxcvbnmQWERTYUIOPASDFGHJKLZXCVBNM.

\multical
QWERTYUIOPASDFGHJKLZXCVBNM.

\multimathop
dist div dom meas sign supp .

\def\accorpa #1#2{\eqref{#1}--\eqref{#2}}
\def\Accorpa #1#2 #3 {\gdef #1{\eqref{#2}--\eqref{#3}}%
  \wlog{}\wlog{\string #1 -> #2 - #3}\wlog{}}

\def\graffe #1{\mathopen\{#1\mathclose\}}

\def\<#1>{\mathopen\langle #1\mathclose\rangle}
\def\norma #1{\mathopen \| #1\mathclose \|}

\def\iot {\int_0^t}

\def\intQt{\int_{Q_t}}
\def\intRt{\int_{R_t}}
\def\intQ{\int_Q}
\def\iO{\int_\Omega}

\def\dt{\partial_t}
\def\dn{\partial_n}

\def\cpto{\,\cdot\,}

\def\checkmmode #1{\relax\ifmmode\hbox{#1}\else{#1}\fi}
\def\aeO{\checkmmode{a.e.\ in~$\Omega$}}
\def\aeQ{\checkmmode{a.e.\ in~$Q$}}

\def\aeS{\checkmmode{a.e.\ on~$\Sigma$}}
\def\aet{\checkmmode{a.e.\ in~$(0,T)$}}

\def\aat{\checkmmode{for a.a.~$t\in(0,T)$}}

\def\erre{{\mathbb{R}}}

\def\genspazio #1#2#3#4#5{#1^{#2}(#5,#4;#3)}
\def\spazio #1#2#3{\genspazio {#1}{#2}{#3}T0}
\def\spaziot #1#2#3{\genspazio {#1}{#2}{#3}t0}

\def\L {\spazio L}
\def\H {\spazio H}
\def\W {\spazio W}
\def\Lt {\spaziot L}

\def\C #1#2{C^{#1}([0,T];#2)}

\def\Lx #1{L^{#1}(\Omega)}
\def\Hx #1{H^{#1}(\Omega)}

\def\Cx #1{C^{#1}(\overline\Omega)}

\def\LQ #1{L^{#1}(Q)}

\def\Ldue{\Lx 2}
\def\Linfty{\Lx\infty}

\def\Huno{\Hx 1}

\def\LQ #1{L^{#1}(Q)}

\let\theta\vartheta
\let\eps\varepsilon
\let\phi\varphi

\let\TeXchi\chi                         % new \chi, exactly on the baseline
\newbox\chibox
\setbox0 \hbox{\mathsurround0pt $\TeXchi$}
\setbox\chibox \hbox{\raise\dp0 \box 0 }
\def\chi{\copy\chibox}

\def\thetaz{\theta_0}
\def\phiz{\phi_0}
\def\thetaQ{\theta_{\!Q}}

\def\Vp{V'}

\def\umin{u_{\rm min}}
\def\umax{u_{\rm max}}
\def\phimin{\phi_\bullet}
\def\phimax{\phi^\bullet}

\def\Uad{\calU_{ad}}
\def\uopt{u^*}
\def\thetaopt{\theta^*}
\def\phiopt{\phi^*}
\def\ubar{\overline u}
\def\thetabar{\overline\theta}
\def\phibar{\overline\phi}
\def\mubar{\overline\mu}
\def\thetah{\theta^h}
\def\phih{\phi^h}
\def\zetah{\zeta^h}
\def\etah{\eta^h}
\def\eh{e^h}

\def\ustar{\bar{u}}
\def\vstar{\bar{v}}

\def\un{u_n}
\def\phin{\phi_n}
\def\thetan{\theta_n}

\def\normaV #1{\norma{#1}_V}
\def\normaH #1{\norma{#1}_H}

\def\normaVp #1{\norma{#1}_*}

\let\hat\widehat
\def\Beta{\hat{\vphantom t\smash\beta\mskip2mu}\mskip-1mu}
\def\betaeps{\beta_\eps}

\def\phieps{\phi_\eps}
\def\xieps{\xi_\eps}

\def\Pi{\hat\pi}

\def\mz{m_0}

\def\xieps{\xi_\eps}

\def\cO{M_\Omega}
\begin{document}
%%%%%%%%%%%%%%%%%%%%%%%%%%%%%%%%%

%%%%%%%%%%%%%%%%%%%%%%%%%%%%%%%%%
%% front page
%%%%%%%%%%%%%%%%%%%%%%%%%%%%%%%%%
%\thispagestyle{empty}

\title{Optimal control for a conserved phase field system\\[0.3cm]
  with a possibly singular potential
}

\author{
Pierluigi Colli\footnote{Dipartimento di Matematica ``F. \pier{Casorati'',} Universit\`{a}
di Pavia, and IMATI-C.N.R., Via Ferrata 5, I-27100 Pavia, Italy (\tt pierluigi.colli@unipv.it).} , 
Gianni Gilardi\footnote{Dipartimento di Matematica \betti{``F. Casorati'',} Universit\`{a}
di Pavia, and IMATI-C.N.R., Via Ferrata 5, I-27100 Pavia, Italy (\tt gianni.gilardi@unipv.it).} ,\\[0.2cm]
Gabriela Marinoschi\footnote{``Gheorghe Mihoc-Caius Iacob'' Institute of Mathematical Statistics and Applied Mathematics of the Romanian Academy,
Calea 13 Septembrie 13, 050711 Bucharest, Romania \gabri{and Simion Stoilow Institute of Mathematics of the Romanian Academy,
Research Group of the Project PN-III-P4-ID-PCE-2016-0372} (\tt gabriela.marinoschi@acad.ro).} \ and\
Elisabetta Rocca\footnote{Dipartimento di Matematica \betti{``F. Casorati'',} Universit\`{a}
di Pavia, and IMATI-C.N.R., Via Ferrata 5, I-27100 Pavia, Italy  
(\tt elisabetta.rocca@unipv.it).}}

\date{}

\maketitle

\begin{abstract}
In this paper we study a distributed control problem for a phase-field system of conserved type with a
possibly singular potential. We mainly handle two cases: the case of a viscous Cahn\pier{--}Hilliard type dynamics for the phase variable in case 
of a logarithmic-type potential with bounded domain and the case of a standard Cahn\pier{--}Hilliard equation in case of a 
regular potential with unbounded domain, 
like the \gianni{classical} double-well potential, for example. Necessary first order conditions of optimality are derived under 
natural assumptions on the data. 
\end{abstract}
\vskip3mm


\pagestyle{myheadings}
\newcommand\testopari{\sc Colli \ --- \ Gilardi \ --- \ Marinoschi \ --- \ Rocca}
\newcommand\testodispari{\sc Optimal control for a conserved phase field system}
\markboth{\testopari}{\testodispari}

%%%%%%%%%%%%%%%%%%%%%%%%%%%%%%%%%
%% very beginning
%%%%%%%%%%%%%%%%%%%%%%%%%%%%%%%%%
%
%  INSERITO \vfill ALTRIMENTI VENIVA SOLO 
%  IL TITOLO DELLA SEZIONE NELLA PRIMA PAGINA 
%
%
%
%

\section{Introduction}
\label{Intro}
\setcounter{equation}{0}

The present contribution is concerned with the study of a distributed control problem for a 
conserved phase field  type PDE system (cf.~\cite{Cag} and \cite{CH58}) in $Q_T := (0,T)\times\Omega $
\Beq
  \dt\theta +\ell \dt\phi - \Delta\theta = u
  \aand
  \dt\phi - \Delta\mu =0, \quad \mu=\tau\dt \phi-\Delta\phi + \calW'(\phi) -\gamma\theta
  \label{caginalp}
\Eeq
where $\Omega$ is the domain where the evolution takes place,
$T$~is some final time,
$\theta$~denotes the relative temperature around some critical value
that is taken to be $0$ without loss of generality,
and $\phi$ is the order parameter.
Moreover, $\ell$ and $\gamma$~are positive coefficients proportional
to the latent heat,
and $u$ is some source term, playing the role of the distributed control here. The parameter $\tau\in [ 0,1]$ 
denotes a viscosity coefficient that  will be taken to be strictly positive or non-negative in the 
subsequent analysis in view of different results. 
Finally, $\calW'$ represents the derivative of a double-well potential~$\calW$,
and the typical example is
the classical regular potential $\calW_{reg}$ defined~by
\Beq
  \calW_{reg}(r) = \frac 14 \, (r^2-1)^2 \,,
  \quad r \in \erre .
  \label{regpot}
\Eeq
However, different choices of $\calW$ are possible,
and a thermodynamically significant example is given by the so-called logarithmic double-well potential,
namely
\Beq
  \calW_{log}(r) = (1+r)\ln (1+r)+(1-r)\ln (1-r) - c r^2 \,,
  \quad r \in (-1,1)
  \label{logpot}
\Eeq
where $c>0$ is large enough in order to kill convexity. 
More generally, the potential $\calW$ could be just the sum $\calW=\Beta+\Pi$,
where $\Beta$ is a convex function that is \gianni{now} allowed to take the value~$+\infty$ in our case 
and $\Pi$ is a smooth perturbation (not necessarily concave).

The mathematical literature on the well-posedness of the PDE system  \eqref{caginalp} is quite vast and so 
we quote here only the papers \cite{BroHilNC}, \cite{CGW, Mir, MZ}, and \cite{KenmNiez} 
dealing respectively with the cases of regular, singular, and non-smooth 
potentials and also with the long-time behavior of solutions.

Moreover, initial conditions like $\theta(0)=\thetaz$ and $\phi(0)=\phiz$ 
and suitable boundary conditions must complement the above equations.
As far as the latter are concerned, 
we take for simplicity the homogeneous Neumann boundary conditions, respectively, that are
\Beq\label{bou}
 \dn  \theta =\dn\phi = \dn \mu=0 
  \quad \hbox{on $\Sigma_T := (0,T)\times\Gamma $} 
\Eeq
where $\Gamma$ is the boundary of~$\Omega$ and $\dn$ is the (say, outward) normal derivative.
We note that the last two boundary conditions are very common in the literature
and that the first one could be replaced by an inhomogeneous~one, for example.
Let us note that by using the third boundary condition in \eqref{bou} we obtain a classical feature of the 
Cahn\pier{--}Hilliard equations, that is the so-called mass conservation:
\[
\int_\Omega \phi(t)=\int_\Omega \phi(0)\quad \forall t\in [0,T]\,. 
\]
The aim of this paper is to study a related optimal control problem for the system \eqref{caginalp}, \eqref{bou}, 
the control being associated to the forcing term $u$ 
that appears on the \rhs\ of the first equation~\eqref{caginalp},
and it  is supposed to vary in some control box~$\Uad$.
We would like to force the averaged temperature and phase variable to be closed to some fixed values $\theta_Q$ and $\phi_Q$ and their final values at time $T$ to be closed to $\theta_\Omega$ and $\phi_\Omega$, respectively. 
In order to do that we choose the following cost functional 
\Beq
  \calJ(u)
  := \frac{\kappa_1}{2}  \intQ (\theta - \thetaQ)^2+ \frac{\kappa_2}{2}  \intQ (\phi - \phi_Q)^2
 + \frac{\kappa_3}{2}  \int_\Omega (\theta(T) - \theta_\Omega)^2+ \frac{\kappa_4}{2}  \int_\Omega (\phi(T) - \phi_\Omega)^2
  \label{Icost}
\Eeq
where $(\theta,\phi)$ is the state corresponding to the control~$u$, and the desired temperatures $\thetaQ\in\LQ2$, $\theta_\Omega\in L^2(\Omega)$, the target 
phases $ \phi_Q \in\LQ2$, $\phi_\Omega\in L^2(\Omega)$, 
and the constants $\kappa_i\geq0$, $i=1,\dots, 4$, are given.
In this case, the optimal control (if~it exists)
balances the smallness of the various differences 
depending on the value of the coefficients~$\kappa_i$.

Thus, the control problem we address in this paper consists in minimizing
the cost functional $\cal J$
depending on the state variables $\theta$ and $\phi$, which satisfy the above state system, 
over all the controls belonging to \gianni{the control box}
\Beq
  \Uad := 
  \bigl\{ u \in L^\infty(Q) : \ \umin\leq u\leq\umax\ \aeQ \bigr\}
  \label{Iuad}
\Eeq
where $\umin$ and $\umax$ are given bounded functions.

The main novelty of the present contribution consists in the fact that we can deal with quite general potentials $\calW$ (even singular) in the phase equation and with a quite general cost functional $\calJ$. Up to our knowledge, indeed, the literature on optimal control for Caginalp type phase field  models is quite poor and often restricted to the case of regular potentials, or dealing with approximating problems when first order optimality conditions are discussed.  In this framework, let us quote the papers \cite{HoffJiang, HKKY} and references therein, as well as\pier{\cite{BBCG, BCF, CGM, CGMR1, CGMR2, CGPS, CGS, CMR, LK, SY, SprZheng}} for different types of phase field  models.
Moreover, up to our knowledge, no optimal control analysis has been performed yet in the literature in case of conserved Capinalp type systems. 
However we can quote the recent results \cite{CFGS, CGS15, CGS16} handling single Cahn\pier{--}Hilliard type dynamics with different boundary conditions and also singular or non-smooth potentials.  

The paper is organized as follows.
In the next section, we list our assumptions, state the problem in a precise form
and present our results. In Sections~\ref{STATE} and~\ref{OPTIMUM}, respectively, we show the well-posedness and regularity results of the state and linearized systems
and the existence of an optimal control. 
The rest (and main part) of the paper is devoted 
to the derivation of first order necessary conditions for optimality.

%%%%%%%%%%%%%%%%%%%%%%%%%%%%%%%%%%%%%%%%%%%%%%%%%%%%%%%%%%%%%%%%%%%%%%%%

\section{Statement of the problem and results}
\label{STATEMENT}
\setcounter{equation}{0}

In this section, we describe the problem under investigation
and present our results. 
As in the Introduction,
$\Omega$~is the body where the evolution takes place.
We assume $\Omega\subset\erre^3$
to~be open, bounded, connected, of class $C^{1,1}$,
and we write $|\Omega|$ for its Lebesgue measure.
Moreover, $\Gamma$ and $\dn$ still stand for
the boundary of~$\Omega$ and the outward normal derivative, respectively.
Given a finite final time~$T>0$,
we set for convenience
\Bsist
  && Q_t := (0,t) \times \Omega 
  \aand
  \Sigma_t := (0,t) \times \Gamma 
  \quad \hbox{for every $t\in(0,T]$}
  \label{defQtSt}
  \\
  && Q := Q_T \,,
  \aand
  \Sigma := \Sigma_T \,.
  \label{defQS}
\Esist
\Accorpa\defQeS defQtSt defQS
Now, we specify the assumptions on the structure of our system.
We assume that
\begin{gather}
  \Beta : \erre \to [0,+\infty]
  \quad \hbox{is convex and lower semicontinuous function with }
  \Beta(0) = 0, 
  \label{hpBeta}
  \\
  \Pi: \erre \to \erre
  \quad \hbox{is a $C^3$ function and 
  ${\Pi\,}'$ is \Lip\ continuous}
  \label{hpPi}
\end{gather}
\Accorpa\HPstruttura hplm hpPi
\gianni{%
and observe that \eqref{hpPi} implies that
\Beq
  |\Pi(r)| \leq \hat c \, (r^2+1)
  \quad \hbox{for every $r\in\erre$}
  \label{disugPi}
\Eeq
with a precise constant~$\hat c$}.
We set for convenience
\Beq
  \gianni{\calW := \Beta + \Pi , \quad}
  \beta := {\Beta}'
  \aand
  \pi := {\Pi\,}'
  \label{defbetapi}
\Eeq
and denote by $D(\beta)$ and $D(\Beta)$ 
the domains of $\beta$ and~$\Beta$, respectively.
We assume then that 
\Beq
  \hbox{$D(\beta)$ is an open interval and $\beta_{|_{D(\beta)}}$ is a $C^2$ function.}
  \label{betareg}
\Eeq
\Accorpa\HPstruttura hpBeta betareg
%Finally, according to the notation in the  Introduction, we set  $\calW:=\Beta+\Pi$. 

We remark that both the regular potential \eqref{regpot} and the logarithmic potential \eqref{logpot}
satisfy the above assumptions on $\beta$ and~$\pi$.
Another possible choice of $\beta$ is given by
\Beq
  \beta(r) := 1 - \frac 1{{r+1}}
  \quad \hbox{for $r>{{}-1}$}
  \label{altrobeta}
\Eeq
and it corresponds to the function $\Beta$ defined by
\Beq
  \Beta(r) := r - \ln (r+1)
  \quad \hbox{if $r>-1$}
  \aand
  \Beta(r) := +\infty
  \quad \hbox{otherwise}
  \label{altroBeta}
\Eeq
with $\Beta$ taking the minimum $0$ at $0$, as required by assumption~\eqref{hpBeta}.
Such an operator $\beta$ yields an example of a different \bhv\ for negative and positive values, 
singular near $-1$ and with a somehow linear growth at $+\infty$.

Moreover, if $\betaeps$~denotes the Yosida regularization of $\beta$ at level~$\eps$,
it is well known that
both $\beta$ and $\betaeps$ are maximal monotone operators
and that $\betaeps$ is even \Lip\ continuous in \gianni{the whole} of~$\erre$. 
Furthermore (see, e.g., \cite[Prop.~2.6, p.~28]{Brezis}), \gianni{we~have}
\Beq
  |\betaeps(r)| \leq |\beta(r)|
  \aand
  \betaeps(r)\to\beta(r)
  \quad \hbox{for $r\in D(\beta)$}.
  \label{propYosida}
\Eeq

Next, in order to simplify notations, we~set
\Beq
  V := \Huno, \quad
  H := \Ldue, \quad
  W := \graffe{v\in\Hx2: \dn v=0}
  \label{defspazi}
\Eeq
and endow these spaces with their natural norms. We have the dense and continuous embeddings 
$W\subset V\subset H\cong H'\subset V'\subset W'$. 
% \betti{We use the same symbol $H$ also to indicate the space $(\Ldue)^3$.} 
We denote by $\langle \cdot, \cdot\rangle_{X', X}$ the duality pairing between two Banach spaces $X'$ and $X$, by $(\cdot, \cdot)_Y$ the scalar product in a generic Hilbert space $Y$, and by $(\cdot,\cdot)$ the scalar product in $H$. Then,  we have $\langle u,v\rangle_{V',V}=(u,v)$ and
$\langle u,w\rangle_{W',W}=(u,w)$ for all $u\in H$, $v\in V$, and $w\in W$. 
The symbol $\norma\cpto_X$ stands for the norm \pier{in a generic Banach space~$X$ or in power of it},
while $\norma\cpto_p$ is the usual norm in both
$\Lx p$ and $\LQ p$, for $1\leq p\leq\infty$. %Moreover, $\norma\cpto_*$ denotes the norm in $V'$.
Finally, for $v\in\L2X$ the function $1*v$ is defined by 
\Beq
  (1*v)(t) := \iot v(s) \, ds
  \quad \hbox{for $t\in[0,T]$}
  \label{defconv}
\Eeq
(note that the symbol $*$ is usually employed for convolution products).

Secondly, we introduce a \wk\ tool, which is useful to deal with a Cahn\pier{--}Hilliard type equation
(see, e.g., \cite[Sect.~2]{CGLN}).
We define the operator
\Beq
  A : V \to \Vp
  \quad \hbox{by} \quad
  \<Av,z>_{V',V} = \iO \nabla v \cdot \nabla z
  \quad \hbox{for every $v,z\in V$}
  \label{defA}
\Eeq
and \gianni{set}
\Beq
  v_\Omega := \frac 1 {|\Omega|} \, \< v , 1 >_{V',V}
  \qquad \hbox{for every $v\in\Vp$}.
  \label{media}
  \Eeq
Recalling our assumption on~$\Omega$, namely, boundedness, smoothness, and connectedness,
we see that the restriction of $A$ to the set of functions $v\in V$ satisfying $v_\Omega=0$
(see~\eqref{media})
is one-to-one and that
$\vstar\in V'$ belongs to the range of $A$ if and only if $\vstar_\Omega=0$.
Therefore, we can define
\Beq
  \dom\calN := \graffe{\vstar\in\Vp: \ \vstar_\Omega = 0}
  \aand
  \calN : \dom\calN \to \graffe{v \in V : \ v_\Omega = 0}
  \label{defN}
\Eeq
by setting:
for $\vstar\in\dom\calN$ and $v\in V$ with $v_\Omega=0$,
the equality $v=\calN\vstar$ means $Av=\vstar$,
i.e., $\calN\vstar$ is the solution $v$ to the \generaliz ed Neumann problem for $-\Delta$
with datum~$\vstar$ that satisfies~$v_\Omega=0$.
This yields a well-defined isomorphism, and the following relations~hold
\begin{align}
  & \iO \nabla \calN \vstar \cdot \nabla v = \< \vstar , v >_{V',V}
  \quad \hbox{for $\vstar\in\Vp$ with $\vstar_\Omega=0$ and $v\in V$}
  \label{dadefN}
  \\
  & \< \ustar , \calN \vstar >_{V',V}
  = \< \vstar , \calN \ustar >_{V',V}
  = \iO (\nabla\calN\ustar) \cdot (\nabla\calN\vstar)\quad
 \hbox{\betti{for $\ustar,\vstar\in\Vp$ with $\ustar_\Omega=\vstar_\Omega=0$}}
  \label{simmN}
  \\[0.2cm]
  & \frac 1 \cO \, \norma{\vstar}_{V'}^2 \leq 
   \normaVp{\vstar}^2:=  \< \vstar , \calN \vstar >_{V',V} \leq \cO \norma{\vstar}_{V'}^2
  \quad \hbox{for all $\vstar\in\Vp$ with $\vstar_\Omega=0$}  \qquad
  \label{normaVp}
\end{align}
for some constant $\cO\geq 1$,
\gianni{%
whence also 
\Beq
  |\< \vstar , v >|
  \leq \cO^{1/2} \normaVp\vstar \normaV v
  \quad \hbox{for all $\vstar\in\Vp$ with $\vstar_\Omega=0$ and $v\in V$.}
  \label{disugdualita}
\Eeq
}%
The first inequality in  \eqref{normaVp} is related to the following Poincar\'e inequality
\Beq
 \normaV{ v }^2\leq \cO (\normaH{\nabla v} + |v_\Omega|)^2
  \quad \hbox{for every $v\in V$}
  \label{poincare}
\Eeq
while \eqref{simmN} implies that we~have
\Beq
  \frac d{dt} \normaVp{\vstar(t)}^2
  = 2 \< \dt\vstar(t) , \calN\vstar(t) >_{V',V}
  \quad \aat
  \label{dtcalN}
\Eeq
for every $\vstar\in\H1\Vp$ satisfying $\vstar_\Omega(t)=0$ for $t\in(0,T)$.

At this point, in order to get useful results both for the state system 
and the linearized one, that we will need later for the optimal control analysis, 
we introduce the following (more general) PDE system which contains the state system 
as particular case. 

Given $\thetaz$ and $\phiz$ such that 
\Bsist
  & \thetaz \in H, \quad \tau^{1/2}\thetaz\in V
    \label{hpz}&
  \\
  & 
  \phiz \in V, \quad
  \Beta(\phiz)\in L^1(\Omega), \quad  m_0:= (\phiz)_\Omega\in D(\beta) &
  \label{hpzbis}
\Esist
and 
\begin{equation}\label{hpvlambda}
v\in L^2(Q), \quad \lambda\in \pier{\H1H \cap L^\infty(Q)} ,
\end{equation}
\Accorpa\HPdati hpz hpvlambda
\Accorpa\HPiniz hpz hpzbis
we look for a triplet $(\theta,\phi,\mu)$ satisfying
\begin{gather}
   \theta \in H^1(0,T;V')\cap \L\infty H\cap \L2 V 
  \label{regtheta}\\
   \tau^{1/2}\theta \in \H1H \cap \L\infty V \cap \L2 W
  \label{regthetabis}
  \\
   \phi \in  H^1(0,T;V')\cap \L\infty V \cap \L2W, \quad \tau^{1/2}\phi\in \H1H
  \label{regphi}
  \\
   \mu\in \L2V,\quad \tau^{1/2}\mu\in\L2W
  \label{regmu}
  \\
 \<\dt\theta+ \ell \dt\phi, z>_{V',V}   +\<A\theta, z>_{V',V}= (v, z) 
  \quad \hbox{ $\forall z\in V$, \, \pier{a.e.~in $(0,T)$}}
  \label{prima}
  \\
   \<\dt\phi, z>_{V',V} +\<A\mu, z>_{V',V} =0 \quad 
   \hbox{ $\forall z\in V$, \, \pier{a.e.~in $(0,T)$}}
   \label{seconda}\\
  \mu= \tau \dt \phi -\Delta\phi + \beta(\phi)+ \gianni\lambda \, \pi(\phi) -\gamma \theta\quad  \aeQ
 \label{terza}
  \\
  % \theta = 0 
  %\aand
  %\dn\phi = 0
  % \quad \aeS
  %\label{bc}
  %\\
   \theta(0) = \thetaz
  \aand
  \phi(0) = \phiz\quad \aeO
  \label{cauchy}
\end{gather}
where the abstract operator $A$ is defined by \eqref{defA}. Note that the initial conditions \eqref{cauchy} make sense since \eqref{regtheta} and \eqref{regphi} entail that $\theta, \, \phi\in \C0H$. 
We also point out that the boundary condition for $\phi$ is included in \eqref{regphi} (cf.~\eqref{defspazi} as well),
\gianni{while those for $\theta$ and $\phi$ are contained in equations \accorpa{prima}{seconda}
due to the definition \eqref{defA} of~$A$}.  
Finally, let us underline that \eqref{seconda}, \eqref{cauchy} and \eqref{hpzbis} easily yield 
\begin{equation}
(\dt \phi)_\Omega =0, \quad \phi_\Omega= m_0 \quad \aet.
\label{pier0}
\end{equation}
\Accorpa\Regsoluz regtheta regmu
\Accorpa\Pblae prima cauchy
\Accorpa\Tuttopblae regtheta cauchy

Our first result, whose proof is sketched in Section~\ref{STATE},
ensures well-posedness with the prescribed regularity, stability
and continuous dependence  in suitable topologies.

\Bthm
\label{Wellposednesslambda}
Assume \HPstruttura\ and \HPdati.
Then, the problem \Pblae\ 
has a unique solution $(\theta,\phi,\mu)$
satisfying \Regsoluz\
and the estimate
\begin{align}
  &\norma\theta_{H^1(0,T;V') \cap \L\infty H\cap \L2V}+\tau^{1/2}\norma\theta_{\H1H \cap \L\infty V \cap \L2W}
  \non
  \\
  &
   +\norma\phi_{ H^1(0,T;V') \cap \L\infty V \cap \L2W}+ \tau^{1/2}\norma\phi_{\H1H}\non\\
&
  + \norma\mu _{\L2V}+\tau^{1/2}\|\mu\|_{\L2W}
  \ \leq \ C_1
  \label{stimasoluz}
\end{align}
holds true for some constant $C_1$ that depends only on $\Omega$, $T$, 
the structure \HPstruttura\ of the system, $\norma{\lambda}_{\pier{\H1H \cap L^\infty(Q)}}$, 
the norms of the initial data associated to~\HPiniz\ and $\norma v_2$.
Moreover, if $v_i\in\LQ2$, $i=1,2$, are given
and $(\theta_i,\phi_i,\mu_i)$ are the corresponding solutions,
then the continuous dependence estimate holds true 
\Bsist
  && \norma{\theta_1-\theta_2}_{\L2H}
  + \norma{(1*\theta_1)-(1*\theta_2)}_{\L\infty V}
  \non
  \\[1mm]
  && \quad {}
  + \norma{\phi_1-\phi_2}_{C^0([0,T]; V')\cap\L2V}+\tau\norma{\phi_1-\phi_2}_{C^0([0,T];H)}
  \non
  \\[1mm]
  && \leq C' \, \norma{(1*v_1)-(1*v_2)}_{\L2H}
  \leq C'' \, \norma{v_1-v_2}_{\L2H}
  \label{contdep}
\Esist
with constants $C'$ and $C''$ that depend only on \pier{$\ell$, $\gamma$,}
$\Omega$, $T$, \pier{$\norma{\lambda}_{L^\infty(Q)}$, and $\norma{\pi'}_{L^\infty (\erre)}$}.
\Ethm

Some further regularity of the solution is stated in the next result, whose proof is given in Section~\ref{STATE}.

\Bthm
\label{Regularitylambda}
The following properties hold true.

$i)$~Assume \HPstruttura\ and \HPdati. Moreover, let \pier{$  v\in L^\infty(Q)$}
\begin{gather}
\label{r01} 
 \pier{\phiz \in W, \quad 
  \beta(\phiz) \in H,\quad  - \Delta \phiz +  \beta(\phiz)  + \lambda(0)\pi(\phiz)\in V}
\\
\label{r02} 
  \pier{\thetaz\in V\cap L^\infty(\Omega)}\, .
\end{gather}
Then, the unique solution $(\theta,\phi,\mu)$ given 
by Theorem~\ref{Wellposednesslambda} also satisfies
\Bsist
&\theta \in \H1H\cap \L\infty V\cap\L2W\cap\LQ\infty&\label{r1}
\\
  & \phi \in W^{1,\infty}(0,T;V') \cap \H1 V \cap \L\infty W, \quad \tau^{1/2}\phi\in W^{1,\infty}(0,T;H)&
  \label{r2}
  \\
  & \mu\in L^2(0,T;W\cap H^3(\Omega))\cap L^\infty(0,T;V)\,\quad \tau^{1/2}\mu \in L^\infty(0,T;W),&
  \label{r3}
  \Esist
  and the initial value (pointwise) problem 
\Bsist
& \dt\theta  + \ell \dt\phi - \Delta\theta= v \quad \aeQ
  \label{primaae}
  \\
  & \dt\phi - \Delta\mu =0 \quad \aeQ
   \label{secondaae}\\
  &\mu= \tau \dt \phi -\Delta\phi + \beta(\phi)+ \gianni\lambda \, \pi(\phi) -\gamma \theta\quad  \aeQ
 \label{terzaae}
  \\
  %& \theta = 0 
  %\aand
  %\dn\phi = 0
  %& \quad \aeS
  %\label{bc}
  %\\
  & \theta(0)=\thetaz
  \aand
  \phi(0) = \phiz\quad \aeO .
  \label{cauchyae}
\Esist
%\Accorpa\Regsoluz regtheta regmu
\Accorpa\Pbl primaae cauchyae
%\Accorpa\Tuttopbl regtheta cauchyae
Besides, the following estimates hold true 
  \Bsist
  &
  \norma\theta_{ \H1H\cap \L\infty V\cap\L2W\cap\LQ\infty}\leq C_2 &
  \label{tetareg}\\
  &\norma\phi_{\W{1,\infty}{V'} \cap \H1V \cap \L\infty W}+\tau^{1/2}\norma\phi_{W^{1,\infty}(0,T;H)}
  \ \leq \ C_3
  \label{phireg}  
  \\
  &\norma\mu_{L^2(0,T;W\cap H^3(\Omega))\cap L^\infty(0,T;V)}+\tau^{1/2}\|\mu\|_{\L\infty W}
\leq C_4\label{mureg}
\Esist
\Accorpa\RegLambda r1 r3
\Accorpa\EstReg tetareg mureg 
for some constants $C_2$, $C_3$, $C_4$ that depend only on $\Omega$, $T$, 
the structure \HPstruttura\ of the system,
the norms of the initial data, \pier{$\norma{v}_\infty$, $\norma{\lambda}_{\pier{\H1H \cap L^\infty(Q)}}$ and the norms of the data in \eqref{r01}--\eqref{r02}.} 

$ii)$~By {further assuming \gianni{that either
$D(\beta)\equiv\erre$ or $\tau>0$ and $\beta(\phiz)\in\Linfty$}},
we have that $\beta(\phi)\in\LQ\infty$ and
\Beq
  \norma{\beta(\phi)}_{\LQ\infty} \leq C_5
  \label{stimaxi}
\Eeq
with a constant $C_5$ that depends on $C_3$, $C_4$, and 
\gianni{even on $\tau$ and $\norma{\beta(\phiz)}_\infty$ if $\tau>0$}.

$iii)$ Moreover, if $\lambda \equiv 1$, $v_i\in\LQ2$, $i=1,2$, are given \pier{and 
$(\theta_i,\phi_i,\mu_i)$ are the corresponding solutions,} % and  $e_i:=\theta_i+\ell\phi_i$,
then the estimate holds true\pier{%
\begin{align}
  & \norma{\theta_1-\theta_2}_{\C0H\cap\L2V}
  +\norma{\phi_1-\phi_2}_{\C0V}
  \non
  \\[1mm]
  & \quad {}
  + \norma{\dt(\phi_1-\phi_2)}_{L^2(\betti{0,T}; V')}+\tau\norma{\dt(\phi_1-\phi_2)}_{\L2H}
  % \non
  % \\[1mm]
  % & 
  \leq C''' \norma{v_1-v_2}_{\L2H}
  \label{contdepbis}
\end{align}
}%
for some constant $C'''$ that depends only on only on 
\gianni{$\ell$, $\gamma$, $\Omega$, $T$, $C_3$, $C_5$, $\beta$ and~$\pi$}.   
\Ethm

By applying \gianni{Theorem~\ref{Wellposednesslambda} and the points~$i)$ and~$ii)$ of Theorem~\ref{Regularitylambda}}
in case $v=u$ and $\lambda=1$\gianni, 
we \pier{deduce} the following existence, uniqueness and regularity results for the 
state system \eqref{caginalp} coupled with boundary conditions \eqref{bou} and initial conditions \eqref{cauchy}.

\Bcor
\label{State}
The following properties hold true.

$i)$ Assume \HPstruttura\ and \HPdati \ with $v=u$ \pier{and $\lambda =1$}.
Then, the following variational formulation of the Cauchy problem associated to the state system \pier{\eqref{caginalp}, \eqref{bou}:}
\Bsist
 & \<\dt\theta+ \ell \dt\phi, z>_{V',V}   +\<A\theta, z>_{V',V} = (u, z) \quad \hbox{ $\forall z\in V$, \pier{a.e.~in $(0,T)$}}
  \label{primastate}
  \\
  & \<\dt\phi, z>_{V',V}  +\<A\mu, z>_{V',V}  =0 \quad \hbox{ $\forall z\in V$, \pier{a.e.~in $(0,T)$}}
   \label{secondastate}\\
  &\mu= \tau \dt \phi -\Delta\phi + \beta(\phi)+ \pi(\phi) -\gamma \theta\quad  \aeQ
 \label{terzastate}
  \\
  %& \theta = 0 
  %\aand
  %\dn\phi = 0
  %& \quad \aeS
  %\label{bc}
  %\\
  & \theta(0)=\thetaz
  \aand
  \phi(0) = \phiz\quad \aeO 
  \label{cauchystate}
\Esist
\Accorpa\PblState primastate cauchystate
 has a unique solution $(\theta,\phi,\mu)$
satisfying \Regsoluz, 
and the estimate \eqref{stimasoluz}
holding true for some constant $C_1$ that depends only on $\Omega$, $T$, 
the structure \HPstruttura\ of the system,
the norms of the initial data associated to~\HPdati\ and $\norma u_2$.
Moreover, if $u_i\in\LQ2$, $i=1,2$, are given
and $(\theta_i,\phi_i,\mu_i)$ are the corresponding solutions,
then the estimate \eqref{contdep}
holds true with constants $C'$ and $C''$ that depend only on \pier{$\ell$, 
$\gamma$,} $T$ and~$\pi$.

$ii)$ Assume \HPstruttura, \HPdati, \eqref{r01}--\eqref{r02} with $v=u$ \pier{and $\lambda =1$. Then,} the unique solution of point $i)$ also satisfies the regularity properties 
\RegLambda, the pointwise system \pier{\Pbl}, and the estimates \EstReg\  with constants 
\pier{depending on $\Omega$, $T$, the structure \HPstruttura\ of the system,
the norms of the initial data, $\| u \|_\infty$ and the norms of the data in \eqref{r01}-- \eqref{r02}.}

$iii)$ Assume moreover \gianni{that either
$D(\beta)\equiv\erre$ or $\tau>0$ and $\beta(\phiz)\in\Linfty$},
we have that $\beta(\phi)\in\LQ\infty$ and \eqref{stimaxi} is satisfied
with a constant $C_5$ that depends on $C_3$, $C_4$, and 
\gianni{even on $\tau$ and $\norma{\beta(\phiz)}_\infty$ if $\tau>0$}.
\Ecor

The well-posedness result for problem \PblState\ given by Corollary~\ref{State} 
allows us to introduce the control-to-state mapping~$\calS$ 
and to address the corresponding control problem.
We define
\begin{align}
  & \calX := L^\infty(Q)
  \label{defX}
  , \quad \ \pier{\calY := (\C0H\cap \L2V)^2}
  \\
  & \calS : \calX \to \calY ,
  \quad 
  u \mapsto \calS(u) =: (\theta,\phi)
  \quad \hbox{where}
  \non
  \\
  & \quad \hbox{$(\theta,\phi)$ is the pair of the first two components} \non\\
&\quad \hbox{of the unique solution  $(\theta,\phi, \mu) $
  to \Regsoluz, \PblState}.
  \label{defS}
\end{align}
Next, in order to introduce the control box and the cost functional,
we assume that
\Bsist
  & \umin,\umax\in\LQ\infty
  \quad \hbox{satisfy} \quad
  \umin {}\leq{}\umax
  \quad \aeQ 
  \label{hpUad}
  \\
  &  \kappa_i \in [0,+\infty), \,  i=1, \dots 4, \quad \displaystyle\sum_{i=1}^4 \kappa_i>0, 
  \quad
  \thetaQ, \, \phi_Q \in \LQ2, \quad \theta_\Omega, \, \phi_\Omega\in H
  \label{hpJ}
\Esist
\Accorpa\StrutturacCP hpUad hpJ
and define $\Uad$ and $\calJ$
according to the Introduction.
Namely, we~set
\Bsist
  && \Uad := 
  \bigl\{ u \in \calX : \ \umin\leq u\leq\umax\ \aeQ \bigr\}
  \phantom \int
  \label{defUad}
  \\[0.2cm]
  && \calJ := \calF \circ \calS : \calX \to \erre
  \quad \hbox{where} \quad
  \calF : \calY \to \erre
  \quad \hbox{is defined by}
  \non % \label{defJ}
  \\
  && \calF(\theta,\phi)
  := \frac{\kappa_1}{2}  \intQ (\theta - \thetaQ)^2+ \frac{\kappa_2}{2}  \intQ (\phi - \phi_Q)^2\non \\
  &&\qquad\qquad\quad
 + \frac{\kappa_3}{2}  \int_\Omega (\theta(T) - \theta_\Omega)^2+ \frac{\kappa_4}{2}  \int_\Omega (\phi(T) - \phi_\Omega)^2 .
  \label{defI}
\Esist
\Accorpa\Defcontrol defX defI

Here is our first result on the control problem; for the proof we refer to Section~\ref{OPTIMUM}.

\Bthm
\label{Optimum}
Assume \HPstruttura, \HPiniz, \eqref{r01}--\eqref{r02}
and let $\Uad$ and $\calJ$ be defined by \accorpa{defUad}{defI}.
Then, there exists $\uopt\in\Uad$ such~that
\Beq
  \calJ(\uopt)
  \leq \calJ(u)
  \quad \hbox{for every $u\in\Uad$} .
  \label{optimum}
\Eeq
\Ethm

Our next aim is to formulate the first order necessary optimality conditions.
As $\Uad$ is convex, the desired necessary condition for optimality~is
\Beq
  ( D\calJ(\uopt) , u-\uopt )_{L^2(Q)}  \geq 0
  \quad \hbox{for every $u\in\Uad$}
  \label{precnopt}
\Eeq
provided that the derivative $D\calJ(\uopt)$ exists 
at least in the G\^ateaux sense  in $L^2(Q)$.
Then, the natural approach consists in proving that 
$\calS$ is \Frechet\ differentiable at $\uopt$
and applying the chain rule to $\calJ=\calF\circ\calS$.
We can properly tackle this project under further assumptions 
on the nonlinearities~$\beta$ and $\pi$.

Since assumptions \HPstruttura\ force $\beta(r)$ to tend to $\pm\infty$
as $r$ tends to a finite end-point of~$D(\beta)$, if any,
we see that combining the further requirements on the initial data 
with the boundedness \pier{properties of} $\phi$ and $\beta(\phi)$ 
stated by Corollary~\ref{State} % at point  $ii)$
immediately yields the following result. 
\Bcor
\label{Bddaway}
Suppose that all the assumptions of Corollary~\ref{State}, point $iii)$ hold true. 
Then, the component $\phi$ of the solution $(\theta,\phi,\mu)$
also satisfies
\Beq
  \phimin \leq \phi \leq \phimax
  \quad \hbox{in $\overline Q$}
  \label{bddaway}
\Eeq
for some constants $\phimin\,,\phimax\in D(\beta)$
that depend only on $\Omega$, $T$, 
the structure  \HPstruttura\ of the system,
the norms of the initial data associated to~\HPiniz,
\gianni{the norms
$\norma u_\infty$, $\norma\thetaz_\infty$, 
and even on $\tau$ and $\norma{\beta(\phiz)}_\infty$ if $\tau>0$}.
\Ecor

As we shall see in Section~\ref{FRECHET}, 
the computation of the \Frechet\ derivative of $\calS$ 
leads to the linearized problem that we describe at once
and that can be stated starting from a generic element $\ubar\in\calX$.
Let $\ubar\in\calX$ and $h\in\calX$ be given.
We set $(\thetabar,\phibar):=\calS(\ubar)$.
Then the linearized problem consists in finding $(\Theta,\Phi, Z)$ 
satisfying 
\begin{gather}
   \Theta \in \H1{H}\cap \pier{\C0V}\cap \L2W\cap L^\infty(Q)
  \label{regTheta}
  \\
   \Phi \in W^{1,\infty}(0,T;V') \cap \H1V\cap \L\infty W, 
   \quad \tau^{1/2}\Phi\in \W{1,\infty}H
  \label{regPhi}
  \\
  Z\in \L2{W\cap H^3(\Omega)}\cap \L\infty V, \quad \tau^{1/2} Z\in \L\infty W
  \label{regZeta}
\end{gather}
and solving the following problem
\Bsist
  & \dt\Theta + \ell\dt\Phi - \Delta\Theta = h
   \quad \aeQ
  \label{linprima}
  \\
  & \dt\Phi - \Delta Z=0
   \quad \aeQ
  \label{linseconda}
  \\
  &Z=\tau \dt \Phi-\Delta\Phi + \calW''(\phibar) \, \Phi -\gamma\Theta
   \quad \aeQ
  \label{linterza}
  \\
   &\dn \Theta = \dn\Phi = \dn Z= 0
  \quad \aeS
  \label{linbc}
  \\
  & \Theta(0) =\Phi(0) = 0
   \quad \aeO .
  \label{lincauchy}
\Esist
\Accorpa\Reglin regTheta regZeta
\Accorpa\Linpbl linprima lincauchy

Applying Theorem~\ref{Regularitylambda} in the case 
$v=h$, $\lambda=\calW''(\phibar)$, $\beta(\phi)=0$, $\pi(\phi)=\phi$, 
\gianni{$\theta_0=0$ and $\phi_0=0$}, 
we \pier{deduce} the 
existence, uniqueness and regularity results for the  linearized system described above.
\gianni{In view of \HPstruttura},
the reader can check that $\calW''(\phibar)$ complies with \eqref{hpvlambda}.

\Bprop
\label{Existlin}
Let the assumptions of Theorem~\ref{Regularitylambda} $ii)$ hold true and let  $\ubar\in\calX$ and $(\thetabar,\phibar)=\calS(\ubar)$.
Then, for every $h\in\calX$,
there exists a unique triplet $(\Theta,\Phi, Z)$
satisfying \Reglin\
and solving the linearized problem \Linpbl.
Moreover, the inequality
\Beq
  \norma{(\Theta,\Phi)}_\calY \leq C_6 \norma h_\calX
  \label{stimaFrechet}
\Eeq
holds true with a constant $C_6$ 
that depend only on $\Omega$, $T$, 
the structure \HPstruttura\ of the system,
the norms of the initial data associated to~\HPiniz,
\gianni{the norms
$\norma\ubar_\infty$, $\norma\thetaz_\infty$,
and even on $\tau$ and $\norma{\beta(\phiz)}_\infty$ if $\tau>0$}.
In particular, the linear map $\calD:h\mapsto(\Theta,\Phi)$
is continuous from $\calX$ to~$\calY$.
\Eprop

In fact, we shall prove that the \Frechet\ derivative 
$D\calS(\ubar)\in\calL(\calX,\calY)$ 
actually exists and coincides with the map $\calD$ introduced in the last statement.
This will be done in Section~\ref{FRECHET}.
Once this is established, we may use the chain rule with $\ubar:=\uopt$ 
to prove that the necessary condition \eqref{precnopt} for optimality takes the form
\Bsist
 && \kappa_1\intQ (\thetaopt-\theta_Q) \Theta+\kappa_2\intQ(\phiopt-\phi_Q)\Phi
 \non\\
 &&+\kappa_3\int_\Omega(\thetaopt(T)-\theta_\Omega)\Theta(T)+\kappa_4\int_\Omega(\phiopt(T)-\phi_\Omega)\Phi(T)
  \geq 0
  \quad \hbox{for any $u\in\Uad$},
  \qquad
  \label{cnopt}
\Esist
where $(\thetaopt,\phiopt)=\calS(\uopt)$ and, for any given $u\in\Uad$, 
the pair $(\Theta,\Phi)$ is the solution to the linearized problem
corresponding to $h=u-\uopt$.

The final step then consists in eliminating the pair $(\Theta,\Phi)$ from~\eqref{cnopt}.
This will be done by introducing  the so-called adjoint problem. 
\Bthm
\label{Existenceadj}
Let the assumptions of Theorem~\ref{Regularitylambda} $ii)$ hold true and let $\uopt$ and $(\thetaopt,\phiopt)=\calS(\uopt)$
be an optimal control and the corresponding state. 
Then there exists a unique solution $(q,p)$ 
with the regularity properties 
\begin{gather}
   q\in \H1{V'} \cap \pier{\C0H} \cap \L2{V}
  \label{regq}
  \\
   p \in \H1{W'} \cap \C0H\cap \L2W, 
  \label{regp}
  \\
  \tau^{1/2} p\in \pier{\H1H \cap  \C0 V}
  \label{regtp} 
\end{gather}
of the adjoint problem
\begin{align}
   -\<\dt q(t), z>_{V',V}
  + \int_\Omega\nabla q(t) \cdot \nabla z
  + \gamma\int_\Omega\Delta p(t) z
  = \int_\Omega \pier{g_1}(t) z \non
  \\
   \forall z\in V, \quad \aat
  \label{primaadj}
\end{align}
\begin{align}
   -\<\dt p(t), w>_{W', W} +\int_\Omega (\tau\dt p(t)+\Delta p(t))\Delta w-\int_\Omega\calW''(\phiopt)\Delta p(t)w +\ell \int_\Omega q(t)\Delta w\non
  \\
  \quad -\ell\gamma\int_\Omega\Delta p(t) w+\int_\Omega(\ell \pier{g_1}(t)-\pier{g_2}(t))w=0 \
  \quad \forall w\in W, \quad \aat
  \label{secondaadj}
\end{align}
\begin{align}
   \<q(T), z>_{V',V} =\int_\Omega \pier{g_3} z\quad \forall z\in V ,\hskip6cm  \non
  \\ 
   \<p(T), w>_{W', W}-\tau \int_\Omega p(T)\Delta w=\int_\Omega (\pier{g_4}-\ell \pier{g_3})w\quad \forall w\in W
  \, 
  \label{cauchyadj}
\end{align}
\Accorpa\Pbladj primaadj cauchyadj
where
\begin{align}
&\pier{g_1}(t)=\kappa_1(\thetaopt(t)-\theta_Q(t)), \quad \pier{g_2}(t)=\kappa_2(\phiopt(t)-\phi_Q(t)), \non
  \\
  &\pier{g_3}=\kappa_3(\thetaopt(T)-\theta_\Omega), \quad \pier{g_4}=\kappa_4(\phiopt(T)-\phi_\Omega).\non
  \end{align}
\Ethm
The proof of the following result will be given in Section~\ref{OPTIMALITY}. 

\Brem
Notice that a strong formulation of \Pbladj \ consists in the following system
\begin{gather}
\gianni{   - \dt q - \Delta q + \gamma\Delta p
  = \kappa_1 (\thetaopt-\thetaQ)
 \quad \aeQ }
  \label{primaadjs}
  \\
\gianni{   - \dt p - \Delta (-\tau\dt p-\Delta p)
  - \calW''(\phiopt)\Delta p - \ell \dt q
  = \kappa_2(\phiopt-\phi_Q)
  \quad\aeQ}
  \label{secondaadjs}
  \\
\gianni{\dn q=\dn p=\dn \Delta p=0
  \quad \aeS}
  \label{bcadjs}
  \\
\gianni{q(T)=\kappa_3(\thetaopt(T)-\theta_\Omega), \quad p(T)-\tau\Delta p(T)+\ell q(T)=\kappa_4(\phiopt(T)-\phi_\Omega)
  \quad \aeO .}
  \label{cauchyadjs}
\end{gather}%
\Accorpa\Pbladjs primaadjs cauchyadjs
\Erem

Our last result, also proved in Section~\ref{OPTIMALITY}, establishes optimality conditions.

\Bthm
\label{CNoptadj}
Let $\uopt$ be an optimal control.
Moreover, let $(\thetaopt,\phiopt)=\calS(\uopt)$ and $(q,p)$
be the associate state and the unique solution to the adjoint problem~\Pbladj\
given by Theorem~\ref{Existenceadj}.
Then we~have
\begin{align}
  \intQ(\uopt-u)q \leq 0 
   \quad \hbox{for every $u\in \Uad$} .   
   \label{cnoptadj}
\end{align}
In particular, we have  $-q\in N_K(\uopt)$, where $K=[u_{\min}, u_{\max}]$ and $N_K$ is the normal cone to the convex set $K$. 
\Ethm

A straightforward consequence of Theorem~\ref{CNoptadj} is here stated.

\Bcor
Under the conditions of Theorem~\ref{CNoptadj}, the
optimal control $\uopt $ reads
\[
u^{\ast }\left\{
\begin{array}{ll}
=u_{\min } &\mbox{a.e. on the set } \ \{(t,x) \, : \ q(t,x)>0\} \\[0.1cm]
=u_{\max } &\mbox{a.e. on the set } \ \{(t,x)\, : \ q(t,x)<0\} \\[0.1cm]
\mbox{$\in(u_{\min}, u_{\max})$}\quad &\mbox{elsewhere.}%
\end{array}%
\right. 
\]
\Ecor

In the remainder of the paper, we often owe to the \Holder\ inequality
and to the elementary Young inequalities
\begin{align}
 ab \leq \alpha \, a^{1/\alpha} + (1-\alpha) \, b^{1/(1-\alpha)}
  \aand
  ab \leq \delta a^2 + \frac 1 {4\delta} \, b^2
  \non
  \\
 \quad \hbox{for every $a,b\geq 0$, \ $\alpha\in(0,1)$ \ and \ $\delta>0$}
  \label{young}
\end{align}
in performing our a priori estimates.
To this regard, in order to avoid a boring notation,
we use the following general rule to denote constants.
The small-case symbol $c$ stands for different constants which depend only
on~$\Omega$, the final time~$T$, the shape of the nonlinearities
and the constants and  norms of
the functions involved in the assumptions of our statements.
A~small-case $c$ with a subscript like $c_\delta$
indicates that the constant might depend on the parameter~$\delta$, in addition.
Hence, the meaning of $c$ and $c_\delta$ might
change from line to line and even in the same chain of equalities or inequalities.
On the contrary, different symbols (e.g., capital letters)
stand for precise constants which we can refer~to.

%%%%%%%%%%%%%%%%%%%%%%%%%%%%%%%%%%%%%%%%%%%%%%%%%%%%%%%%%%%%%%%%%%%%%%%%

\section{The state and the linearized systems}
\label{STATE}
\setcounter{equation}{0}

This section is devoted to the proofs of Theorems~\ref{Wellposednesslambda} and~\ref{Regularitylambda}, which, in turn, imply the 
validity of Corollary~\ref{State} and Proposition~\ref{Existlin}. 
As far as Theorem~\ref{Wellposednesslambda}  is concerned,
we notice that the initial-boundary value problem under study 
is a quite standard phase field system
and that a number of results on it can be found in the literature
(see, e.g., \cite{BroHilNC, BrokSpr, Cag, Gil, Mir}, and references therein).
Nevertheless, we prefer to sketch the basic a~priori estimates
that correspond to the regularity \Regsoluz\ of the solution
and to the stability estimate~\eqref{stimasoluz},
for the reader's convenience.
A~complete existence proof can be obtained by \regulariz ing the problem,
performing similar estimates on the corresponding solution,
and passing to the limit through compactness and monotonicity arguments. 
In particular the potential $\Beta$ should be replaced by 
its \pier{Moreau--}Yosida approximation $\Beta_\eps$, but, since all estimates we deduce 
are formal and independent of $\eps$, we 
skip the index hereby most of the times. 

Concerning the treatment of the unusual term $\lambda (t,x) \pi (\phi)$ in the equation \eqref{terzaae}, we refer the 
reader to the analysis carried out in \cite{CGS1} for a Cahn--Hilliard system with dynamic boundary conditions.

We also give a short proof of~\eqref{contdep} and \eqref{contdepbis}
(whence uniqueness follows as a consequence)
and conclude the discussion on Theorem~\ref{Regularitylambda}.
 
As already mentioned, we derive just formal a priori estimates. Let's define the auxiliary variable $e:=\theta+\ell\phi$. 
We take $z=e$ in  \eqref{prima}; then we test \eqref{seconda} by $L\calN(\dt\phi)$ and \eqref{terza} by $-L\dt\phi$, being $L$ a positive constant 
to be chosen later. Moreover we add to both members of the resulting equality the term $\frac{L}{2}\|\phi(t)\|_H^2+\intQt|\theta|^2$; 
finally, we sum up and integrate over~$Q_t$ with $t\in(0,T)$.
As the terms involving the product $\mu\, \dt\phi$ cancel out, 
 we obtain
\begin{align}
  & \frac 12 \iO |e(t)|^2  
  + \int_0^t\|\theta\|_V^2
  +  L\int_0^t \|\dt\phi\|_{*}^2+\tau L\intQt|\dt\phi|^2+ \frac{L}{2} \norma{\phi (t)}_V^2+ L\iO \Beta(\phi(t))
  \non
  \\
  & = \frac 12 \iO |\thetaz+\ell\phiz|^2
  + \frac{L}{2}\norma{\nabla\phiz}_H^2
   + L\iO \Beta(\phiz)
     + \intQt v \, e -\ell \intQt\nabla\theta \cdot \nabla\phi
     \non\\
     &\quad -L\intQt \gianni\lambda\,\pi(\phi)\,\dt\phi+\gamma L\intQt\theta\,\dt\phi+\frac{L}{2}\|\phi(t)\|_H^2+\intQt|\theta|^2\non\\
     &
=: \frac 12 \iO |\thetaz+\ell\phiz|^2
  + \frac{L}{2}\norma{\nabla\phiz}_H^2
   + L\iO \Beta(\phiz)+ \sum_{i=1}^6 I_i\,.
  \label{s1p}
\end{align}
We can now proceed \pier{by} estimating the six integrals on the right hand side in \eqref{s1p}. 
\pier{Indeed,} the last integral on the \lhs\ is nonnegative thanks to \eqref{hpBeta}
and the first three terms on the \rhs\ are under control, due to \accorpa{hpz}{hpzbis}.
By applying the Young inequality we deduce the estimates 
 \begin{align}
 & I_1\leq \frac12 \int_0^t\|v\|_H^2+\frac12\int_0^t\|e\|_H^2 \label{i1}\\
 &I_2\leq \frac12 \int_0^t\|\nabla\theta\|_H^2+c \int_0^t\|\nabla\phi\|_H^2\,. 
  \label{i2}
\end{align}
\gianni{We treat the third integral by integration by parts in time and taking advantage of the 
continuous embedding $V\subset L^4(\Omega)$. 
Moreover, we account for \gianni{\eqref{disugPi}} and explicitly write the corresponding constant~$\hat c$ in some terms, for clarity.
By allowing the values of $c$ to depend on $L$ as well, we obtain
\Bsist
  && I_3
  = L \intQt \dt\lambda \, \Pi(\phi)
  - L \int_\Omega \lambda(t) \, \Pi(\phi(t))
  + L \int_\Omega \lambda(0) \, \Pi(\phi(0))
  \non
  \\
  && \leq c \intQt |\dt\lambda| \, (|\phi|^2 + 1)
  + L \norma\lambda_\infty \, \hat c \, (\normaH{\phi(t)}^2 + 1)
  + c
  \non
  \\
  && \leq L \norma\lambda_\infty \, \hat c \, \normaH{\phi(t)}^2
  + c \iot \normaH{\dt\lambda} ( \norma\phi_4^2 + 1)
  + c
  \non
  \\
  && \leq L \norma\lambda_\infty \, \hat c \, \normaH{\phi(t)}^2
  + c \iot \normaH{\dt\lambda} \normaV\phi^2 
  + c \,.
  \label{i3}
\Esist
We notice at once that the first summand of the last line is proportional to the term $I_5$
we introduce and treat later~on.
Next, in view of \eqref{disugdualita}}, we have 
\begin{align}
&I_4 = \gamma L\int_0^t \langle \dt\phi, \theta \rangle_{V',V} \leq \frac{L}{4}\int_0^t\|\dt\phi\|_{*}^2+ \gamma^2  L  \cO \int_0^t\|\theta\|_V^2\label{i4}\\
%&I_5= \frac L 2 \norma{\phi_0}^2_H + L \int_0^t  \langle \dt\phi, \phi \rangle_{V',V}
%\leq  \frac{L}{4}\int_0^t\|\dt\phi\|_{*}^2+c \int_0^t\|\phi\|_V^2+c\|\phiz\|_H^2\label{i5}\\
&I_6=  \intQt | e- \ell\phi|^2 \leq c\left(\intQt|e|^2+\intQt|\phi|^2\right)\,.\label{i6}
\end{align}
\gianni{It remains to estimate $I_5:=(L/2)\normaH{\phi(t)}^2$ and the proportional term of~\eqref{i3}.
We observe that
\Beq
  \normaH{\phi(t)}^2
  = \normaH\phiz^2
  + 2 \iot \< \dt\phi , \phi > \,.
  \non
\Eeq
Thus, we have
\Beq
  \bigl( (L/2) + L \norma\lambda_\infty \, \hat c \,\bigr) \normaH{\phi(t)}^2
  \leq \frac L4 \iot \normaVp{\dt\phi}^2
  + c \iot \normaV\phi^2
  + c \,.
  \label{i5}
\Eeq
}%
Choosing now $L$ such that $1- \pier{(1/2) -{}} \gamma^2L\cO  >0$, we insert \gianni{\eqref{i1}--\eqref{i5}} in \eqref{s1p}. Then, 
using \eqref{hpvlambda} together with a standard version of Gronwall lemma, we obtain the following estimate 
\begin{align}
&\|e\|_{\L\infty H}+\|\theta\|_{\L2V}+\|\phi\|_{ \H1 {V'}\cap \L\infty V}\non\\
&+\tau^{1/2}\|\phi\|_{\H1H}+\|\Beta(\phi)\|_{\L\infty {L^1(\Omega)}}\leq c\,.\label{s1}
\end{align}
Hence, by comparison in \eqref{prima} and by virtue of standard regularity results for linear parabolic equations, we have that 
\begin{align}
&\|\dt\theta\|_{\L2{V'}}+\|\theta\|_{\L\infty H}+\tau^{1/2}\|\theta\|_{\H1H\cap \L\infty V\cap \L2W}\leq c\,.
\label{s2}
\end{align}
In view of \eqref{pier0}, we can now test \eqref{seconda} by $\calN(\phi-m_0)$ and subtract \eqref{terza} tested by $\phi-m_0$. 
Two terms cancel out and we can integrate by parts in the term containing~$-\Delta \phi.$ 
\gianni{By rearranging a little, we obtain \aat
\Bsist
  && \iO \betaeps(\phi(t)) (\phi(t)-\mz)
  + \iO |\nabla\phi(t)|^2
  \non
  \\
  && = - \< \dt\phi(t) , \calN(\phi(t)-\mz) >
  - \tau \iO \dt\phi(t) (\phi(t)-\mz)
  \non
  \\
  && \quad {}
  - \iO \lambda(t) \, \pi(\phi(t) \pcol{)}\, (\phi(t)-\mz)
  \pcol{{}+{}} \gamma \iO \theta(t) (\phi(t)-\mz)
  \non
  \\
  && \leq \normaVp{\dt\phi(t)} \, \normaVp{\phi(t)-\mz}
  + \tau \normaH{\dt\phi(t)} \, \normaH{\phi(t)-\mz}
  \non
  \\
  && \quad {}
  + c \norma\lambda_\infty \bigl( \normaH{\phi(t)}^2 + 1 \bigr) 
  + \gamma \normaH{\theta(t)} \, \normaH{\phi(t)-\mz}
  + c \,.
  \label{pier1}
\Esist
Now, we use the fact that $m_0 $ lies in the interior of $D(\beta)$ and consequently (cf.~\cite[Appendix, Prop. A1]{MZ}) 
\Beq
  \betaeps(r) \, (r-\mz) \geq \delta_0 \, |\betaeps(r)| - C
  \non
\Eeq
for every $r\in\erre$ and some positive constants $\delta_0$ and $C$
that do not depend on~$\eps$.}
Hence, thanks to \eqref{s1} we have that
\[
\|\beta_\eps(\phi)\|_{L^2(0,T;L^1(\Omega))}\leq c\,.
\]
Next, by testing \eqref{terza} by $1$, \pier{it is easy to infer that 
\Beq
|\mu_\Omega (t) | \leq \tau \| \dt \phi (t)\|_H + 
\|\beta_\eps(\phi (t) )\|_{L^1(\Omega)} 
 +  c \big(\| \phi (t) \|_H + \| \theta (t)\|_H  + 1 \big)  \label{pier2} 
\Eeq
and} so, by using the estimate (cf.~\eqref{seconda} and \pier{\eqref{normaVp})}
\Beq
\|\nabla (\mu-\mu_\Omega)(t)\|_{2}\leq c\, \norma{\dt \phi (t)}_{V'} 
\label{pier3} 
\Eeq
\pier{for a.a. $t\in (0,T)$,} from \eqref{s1} it follows that 
$$ 
\norma{\mu}_{\L2V} \leq c.
$$
Therefore, we can test \eqref{terza} by $\betaeps (\pier{\phi}) $ and integrate in time; we exploit the nonnegativity of the term $(-\Delta\phi (t) ,  \betaeps (\pier{\phi}(t)) )$\pier{, $\aat$,} in order to recover that 
\[
\|\betaeps(\phi)\|_{\L2H}\leq c
\]
whence, by comparison in \eqref{terza}, we have that $\|\Delta\phi\|_{\L2H} \leq c$. 
From these estimates and by standard elliptic regularity results we \pier{infer the desired} estimate 
\[
\|\phi\|_{\L2W}\leq c\,. 
\]
Let us just comment on the fact that, if we want then to pass to the limit in the 
regularization parameter $\eps$, we can use the strong convergence of the corresponding 
solution $\phi_\eps$ in $\L2V$ which is sufficient, along with the weak convergence of 
$\betaeps(\phieps) $ in $\L2H$, in order \pier{to perform the limit procedure}
in our system. 

Next, we proceed proving  estimate~\eqref{contdep}.
We first integrate \eqref{prima} with respect to time
and get the equation
\Beq
 \<\theta+ \pier{\ell\phi}, z>_{V',V}   +\<A(1*\theta), z>_{V',V}= (\thetaz + \ell \phiz + 1*v, z) \quad \hbox{ $\forall z\in V$, \pier{a.e.~in $(0,T)$}.}
  \label{intprima}
\Eeq
Now, we fix $v_i\in\LQ2$, $i=1,2$,
and consider two corresponding solutions $(\theta_i,\phi_i,\mu_i)$ with the same initial data.
We write \eqref{intprima} for both of them 
and test the difference by \pier{$\gamma \theta/\ell$, where} $\theta:=\theta_1-\theta_2$.
At the same time, we write \eqref{seconda} for both solutions, 
take  the difference and choose $z=\cal N\phi$, where $\phi:=\phi_1-\phi_2$. 
Finally, we take \eqref{terza} 
for the two solutions and test the difference by~\gianni{$-\phi$}. 
Then, we add the resulting equalities and integrate over~$(0,t)$.
\pier{Note that two pairs of corresponding terms cancel.}
Hence, by setting $v:=v_1-v_2$  for brevity, and using the monotonicity of $\beta$, the 
Lipschitz continuity of $\pi$ and the boundedness of $\lambda$, we have
\begin{align}
  & \pier{\frac{\gamma}{\ell}} \intQt |\theta|^2
  + \pier{\frac{\gamma}{2\ell}} \iO |\nabla(1*\theta)(t)|^2 
  +\frac 12 \|\phi(t)\|_{*}^2
  + \frac{\tau}2 \iO | \phi(t)|^2
  +  \intQt |\nabla\phi|^2
  \non
  \\
  &\leq \pier{\frac{\gamma}{\ell}} \intQt (1*v) \, \theta
  \pier{{}- \intQt \gianni\lambda \, \bigl( \pi(\phi_1) - \pi(\phi_2) \bigr) \, \phi } \non
  \\
  & \leq c \norma{1*v}_{\LQ2}^2 + \pier{\frac{\gamma}{2\ell}}\intQt |\theta|^2
  + \norma{\lambda}_\infty \norma{\pi'}_{L^\infty (\erre)} \intQt |\phi|^2 \,.
  \label{pier5}
\end{align}
\pier{Now, we exploit a standard compactness inequality, which states 
that for any $\delta >0$ there is some constant $c_\delta >0$ such that 
\Beq
\|\zeta \|_H^2\leq \delta \|\nabla\zeta \|_{\betti{H}}^2+c_\delta\|\zeta\|_{V'}^2 \,  \quad \hbox{ for all } \, \zeta \in V. 
\label{pier6}
\Eeq
Indeed, by using it to estimate the last term of \eqref{pier5} and owing also to 
\eqref{normaVp}, we have that
$$
\norma{\lambda}_\infty \norma{\pi'}_{L^\infty (\erre)} \intQt |\phi|^2 \leq 
 \frac12 \intQt |\nabla\phi|^2 + 
 c\|\phi(t)\|_{*}^2\betti{.}
$$
Then, by combining it with \eqref{pier5} and applying the} standard Gronwall lemma, 
we obtain the desired estimate \eqref{contdep}.\qed

\medskip

Now, we prove Theorem~\ref{Regularitylambda}.
First take the equation \eqref{prima} and test it by $\dt\theta$, then differentiate \eqref{seconda} and test it by $\calN(\dt\phi)$ and finally take the time derivative of \eqref{terza} and test it by $-\dt\phi$. 
Summing up the resulting equations, \gianni{a~cancellation occurs. So, by} integrating over $(0,t)$, we obtain 
\begin{align}
  &
 \pier{\intQt |\dt\theta|^2} +\frac 12 \|\nabla \theta(t)\|_H^2+ \frac  12 \|\dt\phi(t)\|_*^2+\frac{\tau}{2}\|\dt\phi(t)\|_H^2
  +  \intQt |\nabla\dt\phi|^2
  +  \intQt \beta_\eps'(\phi) |\dt\phi|^2
  \non
  \\
  & \leq \frac 12 \|\nabla \theta_0\|_H^2+
  \frac  12 \|\dt\phi(0)\|_*^2+\frac{\tau}2\|\dt\phi(0)\|_H^2-(\ell-\gamma) \intQt\dt\phi\dt\theta+\intQt v\dt \theta \non
  \\
  &\quad {}
  -\intQt \gianni{\dt\lambda} \, \pi(\phi)\dt\phi-\intQt\lambda(t,x)\pi'(\phi)|\dt\phi|^2
  . 
  \label{s22}
\end{align}
The monotonicity of $\betaeps$ implies that the last term on the \lhs\ is nonnegative\pier{. With the help of \eqref{r01}--\eqref{r02} we find out} that the norms of the 
initial data on the right hand side are bounded\pier{: indeed, write \eqref{seconda}, \eqref{terza} at the time $t=0$, take $z= \calN(\dt\phi(0))$ in \eqref{seconda} and test 
\eqref{terza} by $-\dt\phi (0)$, then sum up and obtain 
$$
\|\dt\phi(0)\|_*^2+\tau\|\dt\phi(0)\|_H^2 \leq \langle \dt\phi (0), \Delta \phiz -  \beta(\phiz)  - \lambda(0)\pi(\phiz) +\gamma \thetaz \rangle_{V',V}
$$
whence
\begin{align}
\pier{\frac 12 \|\dt\phi(0)\|_*^2+\tau\|\dt\phi(0)\|_H^2 \leq c\left(\|\Delta \phiz 
-  \beta(\phiz)  - \lambda(0)\pi(\phiz)\|_V^2+ \|\theta_0\|_V^2\right)}.\non
\end{align}
We can then estimate the next term on the right hand side of \eqref{s22}
by the elementary Young inequality and the 
compactness inequality \eqref{pier6}. Hence, we easily have that}
\[
-(\ell-\gamma) \intQt\dt\phi\dt\theta+\intQt v\dt \theta \leq \frac 12 \intQt|\dt\theta|^2+\delta \intQt|\nabla \dt \phi|^2+c_\delta\int_0^t\|\dt\phi\|_*^2 + \pier{{}c \intQt|v|^2\,}.
\]
\pier{The last two integrals in \eqref{s22} can be treated by means of the regularity assumptions on $\lambda$ and $\pi$ along with the compactness inequality applied to the embedding $V\subset L^4(\Omega)$ as well. We infer that} 
\begin{align}
& -\intQt\lambda_t(t,x)\pi(\phi)\dt\phi-\intQt\lambda(t,x)\pi'(\phi)|\dt\phi|^2\non\\
&\leq \int_0^t\|\lambda_t\|_H\|\pi(\phi)\|_4\|\dt\phi\|_4+c\intQt|\pier{\dt \phi}|^2\non\\
&\pier{\leq \delta \|\nabla \dt \phi\|_{\Lt2H}^2+ c_\delta \|\dt\phi\|_{\Lt2{V'
}}^2+c\int_0^t\|\lambda_t\|_H^2 \, \big( 1+ \|\phi\|_{\L\infty V}^2\big)\,.} \non
\end{align}
Consequently, \pier{taking $\delta$ small enough} we obtain 
\Beq
\norma \theta_{\H1H\cap \L\infty V}+ \norma\phi_{\W{1,\infty}{V'}\cap\H1V}+\tau^{1/2}\|\phi\|_{\W{1,\infty}H}
  \leq c \,.
  \label{pier4}
\Eeq
\pier{At this point, we go back to \eqref{prima} and observe that a comparison of terms entails 
$\norma{A\theta}_{\L2H} \leq c$, whence (cf.~\eqref{defA}) 
\Beq
\norma \theta_{\L2W}
  \leq c 
  \non
\Eeq
and \eqref{primaae} holds.}
Now, \pier{since $\dt\phi$ is bounded in $\L2{L^6(\Omega)}$, $v$ is in $L^\infty(Q)$ and $\thetaz \in L^\infty (\Omega)$, from \eqref{primaae} and the parabolic regularity theory (cf.~\cite[Thm.~7.1, p.~181]{LSU}) it is straightforward to infer that}
\[
\norma\theta_{L^\infty(Q)}\leq c\,. 
\]
\pier{In view of \eqref{pier4} and recalling the estimates \eqref{pier1}--\eqref{pier3} we 
easily conclude that 
$$ 
\norma{\betaeps (\phi) }_{\L\infty {L^1(\Omega)}} + \|\mu\|_{\L\infty V}\leq c\,. 
$$
Next,} by comparison in \eqref{terza} we \pier{obtain} that the term \pier{$-\Delta\phi+\beta_\eps(\phi)$ is bounded in $\L\infty H$, then it is now a standard matter to check that both 
$\|\beta_\eps(\phi)\|_{\L\infty H}$ and 
$\|\Delta\phi\|_{\L\infty H}$
are bounded, whence $\|\phi\|_{\L\infty W} \leq c$ 
on account of \eqref{pier4} as well.}
Moreover, \pier{as we are working in 3D and $W$ is complactly embedded in~$\Cx0$,
from, e.g., \cite[Sect.~8, Cor.~4]{Simon} it follows that 
$\phi$ is bounded in $\C0{\Cx0}=C^0(\overline Q)$.}
Finally, 
we observe that \pier{from \eqref{pier4} and \eqref{seconda} 
it is easy to deduce that
\[
\norma{\mu}_{\L2 {W\cap H^3(\Omega)}} + \norma{\tau^{1/2}\mu}_{\L\infty W}\leq c 
\]
and consequently
$
\norma{\tau^{1/2}\mu}_{L^\infty(Q)}\leq c .
$ This proves~$i)$.}

For the second statement~$ii)$, we can write \pier{\eqref{terza}} in the form
\Beq
  \tau\dt\phi - \Delta\phi + \xi = f:= \pier{\mu +{}} \gamma\theta - \lambda(t,x)\pi(\phi) ,
  \pier{\quad \hbox{with} \quad}
  \xi =\beta(\phi),
  \quad \aeQ
  \label{secondabis}
\Eeq
and \pier{observe} that $\tau^{1/2}f$ is bounded in~$\LQ\infty$
on account of the result $i)$ just proved\pier{. Then,} we 
can use the same estimate already performed in \cite{CGMR1}, 
i.e., we can \pier{multiply} the approximation of \eqref{secondabis} 
\Beq
  \pier{\tau}  \dt\phieps - \Delta\phieps + \xieps = f
  \pier{\quad \hbox{with} \quad}
  \xieps := \betaeps(\phieps), 
  \quad \aeQ
  \label{secondaeps}
\Eeq
\pier{by $|\xieps|^{p-1}\sign\xieps$,} where $\betaeps$ is the Yosida \regulariz ation of~$\beta$ at level~$\eps>0$ \pier{and} $p>2$ is arbitrary, and integrate over~$Q_t$. 
Indeed, a~standard argument shows that $\phieps$ converges to $\phi$ 
in the proper topology as $\eps$ tends to zero,
so that $ii)$~immediately follows
whenever we prove that $\xieps$ is bounded in $\LQ\infty$ uniformly with respect to~$\eps$.
This estimate leads \pier{plainly} to 
\[
\tau^{1/2}\|\beta(\phi)\|_{L^\infty(Q)}\leq c\,. 
\]
Hence, in case $\tau>0$, the proof of \gianni{$ii)$ of} Theorem~\ref{Regularitylambda} is completed. 
In case $\tau=0$ and assuming that $D(\beta)\equiv\erre$, the 
\pier{boundedness of  $\norma{\beta(\phi)}_{L^\infty(Q)}$ 
is an easy consequence of  the facts} that $\phi$ 
is bounded in $C^0(\overline Q)$ and the \pier{real function}
$\beta$ is bounded on bounded sets. 

We need now to prove $iii)$, that is the continuous dependence estimate \eqref{contdepbis},
\gianni{and a preliminary remark is needed.
As pointed out before its statement,
Corollary~\ref{State} depends only on Theorem~\ref{Wellposednesslambda} and on the points~$i)$ and~$ii)$ of Theorem~\ref{Regularitylambda}.
The same holds for Corollary~\ref{Bddaway} as a consequence.
Therefore, in proving~\eqref{contdepbis}, we can use \eqref{bddaway} for every solution.
In particular, we can assume $\calW'$ and $\calW''$ to be \Lip\ continuous and bounded without loss of generality.}
Let's define $v:=v_1-v_2$. We test the \pier{difference} of \eqref{prima} corresponding to different 
solutions $(\theta_i,\phi_i)$\pier{, $i=1,2,$ by $\theta:=\theta_1-\theta_2$, 
the difference} of \eqref{seconda} corresponding to different 
solutions $(\theta_i,\phi_i)$\pier{, $i=1,2,$} by $M\calN(\dt\phi):=M\calN(\dt(\phi_1-\phi_2))$,  the \pier{difference} of \eqref{terza} corresponding to different 
solutions $(\theta_i,\phi_i)$\pier{, $i=1,2,$} by \gianni{$-M\dt\phi:=-M\dt(\phi_1-\phi_2)$}, 
with $M$ \pier{chosen equal to $\ell/\gamma$ in order to cancel two terms in the sum.
We} integrate over $(0,t)$ and sum the three resulting equations up\pier{, thus obtaining
\begin{align}
&\frac 12 \norma{\pier{\theta}(t)}_H^2+\intQt \pier{|\nabla\theta|^2} +M\int_0^t\|\dt\phi\|_*^2+M\tau\int_0^t\|\dt\phi\|_H^2\non\\
&+\frac{M}{2}\|\nabla\phi(t)\|_H^2=\intQt v \, \pier{\theta} -\intQt M(\calW'(\phi_1)-\calW'(\phi_2))\dt\phi
%+\gamma M\intQt\theta\dt\phi=:\sum_{i=1}^3 I_i
\,.
\label{ss3}
\end{align} 
Now, we have that 
$$
\intQt v \, \pier{\theta}\leq \frac 12 \|v\|_{\L2H}^2+\frac 12 \|\theta\|_{\Lt2H}^2\non\\
$$
and the last integral on the right hand side of \eqref{ss3} can be estimated by using 
\eqref{normaVp},  \eqref{phireg} and \eqref{stimaxi} as follows
\gianni{(where the values of $c$ can depend on~$M$)}:  
\begin{align*}
&- \intQt M(\calW'(\phi_1)-\calW'(\phi_2))\dt\phi \non\\
&\leq \frac M 2 \int_0^t\|\dt\phi\|_*^2 + c \|\calW'(\phi_1)-\calW'(\phi_2)\|_{\L2V}^2\non\\
&\leq \frac M 2 \int_0^t\|\dt\phi\|_*^2
\non \\
&\quad+c\left(\|\phi\|_{\L2H}^2+\intQt|(\calW''(\phi_1)-\calW''(\phi_2))\nabla\phi_1|^2+\intQt|\nabla\phi\calW''(\phi_2)|^2\right)\non
\\
&\leq \frac M 2 \int_0^t\|\dt\phi\|_*^2 +c \left(\|\phi\|_{\L2H}^2+\intQt|\phi|^2|\nabla\phi_1|^2+ \betti{\norma{\calW''(\phi_2)}^2_\infty}  \intQ|\nabla\phi|^2\right)\non\\
& \leq \frac M 2 \int_0^t\|\dt\phi\|_*^2+c\left(\|\phi\|_{\L2V}^2+
\|\phi_1\|^2_{\L\infty {W^{1,4}(\Omega)}} 
\int_0^T\|\phi\|_4^2\right)\non\\
& \leq \frac M 2 \int_0^t\|\dt\phi\|_*^2+ c\left(  1 + \|\phi_1\|^2_{\L\infty W}\right)
\|\phi\|_{\L2V}^2.
\end{align*}
Hence, thanks to the already shown estimate \eqref{contdep}, from \eqref{ss3} we infer that 
\begin{align*}
&\frac 12 \norma{\pier{\theta}(t)}_H^2+\intQt \pier{|\nabla\theta|^2} +\frac M 2\int_0^t\|\dt\phi\|_*^2+M\tau\int_0^t\|\dt\phi\|_H^2
+\frac{M}{2}\|\nabla\phi(t)\|_H^2 \\
&\leq 
\frac 12 \|\theta\|_{\Lt2H}^2  +   c \|v\|_{\L2H}^2 .
\end{align*}
Then, by applying the Gronwall lemma we end up with the desired estimate \eqref{contdepbis}. This concludes the proof of Theorem~\ref{Regularitylambda}.} 
\qed

%%%%%%%%%%%%%%%%%%%%%%%%%%%%%%%%%%%%%%%%%%%%%%%%%%%%%%%%%%%%%%%%%%%%%%%%

\section{Existence of an optimal control}
\label{OPTIMUM}
\setcounter{equation}{0}

The following section is devoted to the proof of  Theorem~\ref{Optimum}. We use the direct method,
observing first that $\Uad$ is nonempty.
Then, we let $\graffe{\un}$ be a minimizing sequence for the optimization problem
and, for any~$n$, we take the corresponding solution $(\phin,\thetan,\mu_n)$ to problem~\Pbl.
Then, $\graffe{\un}$ is bounded in \pier{$L^\infty (Q)$} and estimates \pier{\EstReg\ hold} for $(\phin,\thetan,\mu_n)$.
Therefore, we have for a subsequence
\Bsist
  & \un \to u
  & \quad \hbox{weakly star in \pier{$L^\infty (Q)$}}
  \non
  \\
  & \thetan \to \theta
  & \quad \hbox{weakly star in $\H1H \cap \L\infty V \cap L^\infty(Q)$}
  \non
  \\
  & \phin \to \phi
  & \quad \hbox{weakly star in $\W{1,\infty}{V'} \cap \H1V \cap \L\infty W$}
  \non
  \\
  &\tau^{1/2}\phin \to \tau^{1/2}\phi
  & \quad \hbox{weakly star in $\W{1,\infty}{H}$}
  \non
  \\
  & \mu_n\to \mu
  & \quad \hbox{weakly star in $\pier{\L\infty V \cap L^2(0,T;W \cap H^3(\Omega))}$} 
  \non
  \\
  &\tau^{1/2}\mu_n\to \tau^{1/2}\mu & \quad \hbox{weakly star in $\L\infty W$} 
  \non
\Esist
\pier{and $\beta (\phin)$ converges to some $\xi$ weakly star in $\L\infty H$.
Then, in view of \eqref{defUad} it is clear that $u\in\Uad$,}
the initial conditions for $\theta$ and $\phi$ are satisfied,
and we can easily conclude by standard \pier{arguments}.
Very shortly, $\graffe{\phin}$~converges strongly, e.g., in~$L^2(Q)$ and \aeQ\ (for~a subsequence)
by~the Aubin-Lions compactness lemma (see, e.g., \cite[Thm.~5.1, p.~58]{Lions}),
whence $\pi(\phin)$ converges to $\pi(\phi)$ is the same topology
and \pier{$\beta(\phin)\to \xi= \beta(\phi)$ by the weak-strong convergence property} (see, e.g., \cite[Lemma~1.3, p.~42]{Barbu}).
Thus, $(\theta,\phi,\mu)$ satisfies problem~\Pbl.
On the other hand, $\calF(\thetan,\phin)$ converges both to the infimum of $\calJ$
and to $\calF(\theta,\phi)$. 
Therefore, $u$~is an optimal control.\qed

%%%%%%%%%%%%%%%%%%%%%%%%%%%%%%%%%%%%%%%%%%%%%%%%%%%%%%%%%%%%%%%%%%%%%%%%

\section{The control-to-state mapping}
\label{FRECHET}
\setcounter{equation}{0}

As sketched in Section~\ref{STATEMENT}, 
the main point is the \Frechet\ differentiability of the control-to-state mapping~$\calS$.
This involves the \lineariz ed problem \Linpbl,
whose well-posedness is stated in Proposition~\ref{Existlin}.

Here is the main result of this section.

\Bthm
\label{Fdiff}
Let $\ubar\in\calX$ and let \pier{$\calS(\ubar)$ be  
the pair $(\thetabar,\phibar)$ of the first two components 
of the unique solution  $(\thetabar,\phibar, \mubar) $
to \Regsoluz, \PblState\ \gianni{with $u=\ubar$}. 
Then, $\calS$~is \Frechet\ differentiable at~$\ubar$
and the \Frechet\ derivative $[D\calS](\ubar)$}
\pier{is precisely} the map $\calD\in\calL(\calX,\calY)$ defined
in the statement of Proposition~\ref{Existlin}.
\Ethm

\Bdim
We fix $\ubar\in\calX$ and the corresponding state $(\thetabar,\phibar)$
and, for $h\in\calX$ with $\norma h_\calX\leq\Lambda$, for some positive constant $\Lambda$, we~set
\Beq
  (\thetah,\phih) := \calS(\ubar+h)
  \aand
  (\zetah,\etah, \xi^h) := (\thetah-\thetabar-\Theta,\phih-\phibar-\Phi, \mu^h-\mubar-Z)
  \non
\Eeq
where $(\Theta,\Phi, Z)$ is the solution to the linearized problem corresponding to~$h$.
We have to prove that
$\norma{(\zetah,\etah)}_\calY/\norma h_\calX$ tends to zero as $\norma h_\calX$ tends to zero.
More precisely, we show that
\Beq
  \norma{(\zetah,\etah)}_\calY
  \leq c \norma h_{\LQ2}^2
  \label{tesiFrechet}
\Eeq
for some constant~$c$,
and this is even stronger than necessary.
First of all, we fix one fact.
As both $\norma\ubar_\infty$ and $\norma{\ubar+h}_\infty$
are bounded by $\norma\ubar_\infty+\Lambda$,
we can apply Corollary~\ref{Bddaway}
and find constants $\phimin,\phimax\in D(\beta)$ such that
\Beq
  \phimin \leq \phibar \leq \phimax
  \aand
  \phimin \leq \phih \leq \phimax
  \quad \aeQ .
  \label{perstimaTaylor}
\Eeq
Now, let us prove~\eqref{tesiFrechet} by writing the problem solved by~$(\zetah,\etah)$.
We clearly have
\Bsist
  & \dt\zetah - \Delta\zetah + \ell \dt\etah = 0
 \quad \aeQ
  \label{primah}
  \\
  & \dt\etah - \Delta\xi^h=0 \quad \aeQ\label{secondah}\\
  &\xi^h=\tau\dt\etah \gianni{{}-{}} \Delta\etah+\calW'(\phih) - \calW'(\phibar) - \calW''(\phibar) \, \Phi -\gamma \zetah
  \quad \aeQ .
  \label{terzah}
\Esist
Moreover, $\zetah$, $\etah$, and $\xi^h$  satisfy  all homogeneous Neumann boundary conditions and $\zetah$, $\etah$ satisfy homogeneous initial conditions.
At this point, we multiply \eqref{primah}  by $\zetah+\ell\etah$ and sum it up to \eqref{secondah} tested by $\tilde\ell \calN\etah$ and to \eqref{terzah} tested by $-\tilde \ell \etah$, with $\tilde \ell$ a positive constant to be chosen later. 
\gianni{The terms involving $\xi^h$ cancel each other. Thus, integrating}
 the resulting equality over~$(0,t)$, we obtain 
\Bsist
  && \frac 12 \|(\zetah+\ell\etah)(t)\|_H^2+\intQt |\nabla \zetah|^2
  + \frac{\tilde\ell}2 \|\etah(t)\|_*^2+\frac{\tau\tilde\ell}{2}\|\etah(t)\|_H^2
  \non
  \\
  && \pier{+  \tilde\ell\intQt |\nabla\etah|^2} = -\ell\intQt\nabla\zetah\cdot\nabla\etah-  \intQt \tilde \ell I^h\etah-\intQt\gamma\tilde\ell \zetah\etah \,,
  \label{dastimare}
\Esist
where we have defined 
\[
I^h=\calW'(\phih)-\calW'(\phibar)-\calW''(\phibar)\Phi=\calW''(\phibar)\etah+\frac 12 \calW'''(\tilde\phi_h)(\phih-\phibar)^2\,,
\]
\gianni{\pcol{$\tilde\phi_h$ being} some function whose values lie between those of $\phi^h$ and $ \bar \phi$.
In particular, the analogue of \eqref{perstimaTaylor} holds for~$\tilde\phi_h$, so that $\calW'''(\tilde\phi_h)$ is bounded.
The same} is true for~$\calW''(\phibar)$. 
Now we can \gianni{deduce an estimate for} the \rhs\ of \eqref{dastimare}
by accounting for the Young and \Holder\ inequalities,
\gianni{the compactness inequality \eqref{pier6} and the continuous embedding $V\subset L^4(\Omega)$. 
We first observe that 
\Beq
  - \ell \intQt \nabla\zetah \cdot \nabla\etah
  \leq \frac 12 \intQt |\nabla\zetah|^2 + \frac {\ell^2} 2 \intQt |\nabla\etah|^2 .
  \non
\Eeq
Therefore,}
\pier{letting $\tilde\ell> \ell^2 / 2$, 
setting $L=\tilde\ell -  \ell^2 / 2$, \gianni{defining $\eh=\zetah+\ell\etah$ and adding the term $L\intQt|\etah|^2$ to both sides,}
we have that}
\begin{align}
  &\frac 12 \|\eh(t)\|_H^2+\frac 12 \intQt|\nabla\zetah|^2+\frac{\tilde\ell}{2}\|\etah(t)\|_*^2+\frac{\tilde\ell\tau}{2}\|\etah(t)\|_H^2+L\int_0^t\|\etah\|_V^2\non\\
  &\leq \gianni- \intQt\tilde \ell I^h\etah-\gianni{\gamma\tilde\ell}\intQt\zetah\etah
  \pier{{}+ L\intQt|\etah|^2}
  \non\\
  &\leq \intQt  \pier{\big( \tilde\ell\, \calW''(\phibar) + L\big)|\etah|^2}+\frac{\tilde\ell}{2}\intQt\calW'''(\tilde\phi_h)(\phih-\phibar)^2\etah+\intQt\gamma\ell\tilde\ell(\etah)^2-\intQt\gamma\tilde\ell\etah\eh
  \non\\
  &\leq  c\intQt|\etah|^2+c\int_0^t \pier{\|\phih-\phibar\|_4^2 \, \|\etah\|_H}  +c\left(\intQt|\etah|^2+\intQt|\eh|^2\right)\non
  \\
  &\leq \frac{L}{2}\int_0^t\|\etah\|_V^2+c\left(\int_0^t\|\etah\|_*^2+\int_0^t\|\eh\|_H^2\right)+\int_0^t\pier{\|\phih-\phibar\|_V^4}\,. \non
\end{align}
Now, we recall that estimate \eqref{contdepbis} holds for the pair of controls $\ubar+h$ and $\ubar$
and for the corresponding states 
$(\thetah,\phih)$ and $(\thetabar,\phibar)$.
Therefore, we can \pier{proceed} and obtain
\[
\pier{\int_0^t\|\phih-\phibar\|_V^4 \leq c\|\phih-\phibar\|_{\L\infty V}^4\leq \|h\|_{L^2(Q)}^4 . }
\]
\pier{Then, the application of the Gronwall lemma} closes the estimate and yields 
\begin{align}
& \pier{\|\eh(t)\|_H^2}+\intQt |\nabla \zetah|^2
  + \|\etah(t)\|_*^2+\tau\|\etah(t)\|_H^2
  +  \pier{\int_0^t\|\etah\|_V^2}
 \leq c\|h\|_{L^2(Q)}^4\, \label{frechet1}
\end{align}
\pier{\aat .} In order to conclude the proof of \eqref{tesiFrechet}, we need an estimate  in $\C0H$ and so we test \eqref{secondah} by $\etah$ and \pier{add it} to \eqref{terzah} tested by $\Delta\etah$. 
Integrating over $(0,t)$  and using Young's inequality with \eqref{frechet1}, 
we obtain 
\[
\frac12 \|\etah(t)\|_H^2 
+\frac \tau 2  \int_\Omega|\nabla\etah (t)|^2
+ \frac12 \intQt |\Delta \etah  |^2
\leq c \intQt |I_h - \gamma (\eh - \ell\etah) |^2 \leq c \|h\|_{L^2(Q)}^4\]
and by comparison, we also get 
\[
\|\zetah(t)\|_H^2\leq c \|h\|_{L^2(Q)}^4 
\]
\aat, \pier{which} concludes the proof
\gianni{since $\norma h_{\LQ2}\leq c \,\norma h_{\calX}$}.
\Edim

\gianni{%
\Brem
\label{SceltaX}
We have choosen $\calX=\LQ\infty$ by~\eqref{defX}.
However, the $L^\infty$ norm has been used just at the beginning of the proof 
and some modification is possible.
In particular, we can make the more suitable choice $\calX=\LQ2$ and perform the same argument to prove
the directional differentiability of $\calS$ in all the directions $h\in\LQ\infty$.
Indeed, $\ubar\in\LQ\infty$ since $\ubar\in\Uad$.
We point out that this modification does not have any bad consequence in the results of the next section, 
since the necessary condition we prove only uses the directional differentiability of~$\calJ=\calF\circ\calS$,
which still holds in the modified framework.
\Erem
}

%%%%%%%%%%%%%%%%%%%%%%%%%%%%%%%%%%%%%%%%%%%%%%%%%%%%%%%%%%%%%%%%%%%%%%%%

\section{Necessary optimality conditions}
\label{OPTIMALITY}
\setcounter{equation}{0}

In this section, we derive the optimality condition~\eqref{cnoptadj}
stated in Theorem~\ref{CNoptadj}.
We start from~\eqref{precnopt} and first prove~\eqref{cnopt}.

\Bprop
\label{CNopt}
Let $\uopt$ be an optimal control and $\pier{(\thetaopt,\phiopt)}:=\calS(\uopt)$.
Then \eqref{cnopt} holds.
\Eprop

\Bdim
This is essentially due to the chain rule for \Frechet\ derivatives, 
as already said in Section~\ref{STATEMENT},
and we just provide some detail.

It follows that $\calF$ is \Frechet\ differentiable in $\calZ:=\pier{\C0H\times\C0H}$
and that its \Frechet\ derivative $[D\calF](\thetabar,\phibar)$ 
at any point $(\thetabar,\phibar)\in\calZ$ acts as follows
\begin{align}
 [D\calF](\thetabar,\phibar):
  (h_1,h_2)\in\calZ \mapsto 
 \quad \ &\kappa_1  \intQ (\theta - \thetaQ)\pier{h_1}+ \kappa_2\intQ (\phi - \phi_Q)\pier{h_2}
 \non \\
 &+ \kappa_3  \int_\Omega (\theta(T) - \theta_\Omega)\pier{h_1(T)}+ \kappa_4 \int_\Omega (\phi(T) - \phi_\Omega)\pier{h_2(T)} \,.
  \non
\end{align}
Therefore, Theorem~\ref{Fdiff} and the chain rule ensure that
$\calJ$ is \Frechet\ differentiable at $\uopt$
and that its \Frechet\ derivative $[D\calJ](\uopt)$ \pier{at} any optimal control $\uopt$
\pier{is specified by}
\begin{align}
  [D\calJ](\uopt):
  h \in \calX \mapsto\ 
   \quad &\kappa_1  \intQ (\theta - \thetaQ)\Theta+ \kappa_2\intQ (\phi - \phi_Q)\Phi
   \non \\
 &+ \kappa_3  \int_\Omega (\theta(T) - \theta_\Omega)\Theta(T)+ \kappa_4 \int_\Omega (\phi(T) - \phi_\Omega)\Phi(T) 
  \non
\end{align}
where $(\Theta,\Phi)$ is the solution to the linearized problem corresponding to~$h$.
Therefore, \eqref{cnopt} immediately follows from~\eqref{precnopt}.
\Edim

\pier{The next step is the proof of Theorem~\ref{Existenceadj}.
As far as existence is concerned, we can  derive a basic formal estimate.
We take as test functions \pier{$z=q$ in \eqref{primaadj}, $w=p$ in \eqref{secondaadj}
and add the equalities we obtain. Then, we integrate over $(t,T)$ using 
the final conditions \eqref{cauchyadj}.} This computation leads to
\begin{align}
&\frac12 \iO |q(t)|^2 + \intRt |\nabla q|^2 
+ \frac12\iO |p(t)|^2 + \frac\tau 2 \iO |\nabla p(t)|^2 + \intRt |\Delta p |^2 
\non \\  
&= \frac12\int_\Omega |\pier{g_3}|^2 + \frac12\int_\Omega |\pier{g_4}-\ell \pier{g_3}|^2
 -(\gamma +\ell)\intRt q\, \Delta p
\non \\  
&\quad + \intQt (\calW''(\phiopt) -\ell\gamma) p \,\Delta p
  +\intRt \pier{g_1} q -\intQt(\ell \pier{g_1}-\pier{g_2})p 
  \label{pier8}
\end{align}
where $R_t := (t,T) \times \Omega$. We observe that $\calW''(\phiopt)$ is uniformly bounded 
in view of Corollary~\ref{Bddaway} and due to the properties \HPstruttura\ of $\beta$.
Hence, recalling the definitions of $g_1, \, \ldots, g_4$ and owing to the Young inequality~\eqref{young}, we easily infer that 
\begin{align}
&\frac12 \iO |q(t)|^2 + \intRt |\nabla q|^2 
+ \frac12\iO |p(t)|^2 + \frac\tau 2 \iO |\nabla p(t)|^2 + \frac12 \intRt |\Delta p |^2 
\non \\  
  &\leq c\left(\intQt|p|^2+\intQt|q|^2+\|\thetaopt\|_{\C0H}^2+\|\phiopt\|_{\C0H}^2\right)\non
  \\
  &\quad+c\left(\|\theta_Q\|_{L^2(Q)}^2+\|\phi_Q\|_{L^2(Q)}^2+\|\theta_\Omega\|_H^2+\|\phi_\Omega\|_H^2\right)\,. 
  \non
\end{align}
Therefore, we can apply the Gronwall lemma and deduce that 
\Beq
  \norma q_{\C0H\cap\L2V}   +  \norma p_{\C0H\cap \L2W } + \tau^{1/2} \norma p_{\C0V } 
  \leq c\,.
  \label{formala}
\Eeq
This procedure implies in particular the uniqueness of the solution, due to the linearity of the problem:
indeed,  we can replace all $\pier{g_i}$'s in \eqref{pier8} by~$0$ for the difference of two solutions.
Moreover, in the light of \eqref{formala} we can compare the terms of 
\eqref{primaadj} and \eqref{secondaadj} and deduce the estimate
\Beq
  \norma{ \dt q}_{\L2{V'}}   +  \norma {\dt p + \tau \Delta \dt p  }_{\L2{W'} } \leq c 
  \label{formala2}
\Eeq
which enables us to recover the full regularity of the solution in \eqref{regq}--\eqref{regtp}.
Therefore, it is clear how to give a rigorous proof
based on a Faedo--Galerkin scheme, by choosing a basis of eigenfuntions related to the 
operator $- \Delta $ with Neumann homogeneous boundary conditions (cf.~\eqref{defA}).  
This approximation scheme would provide 
a sequence $\graffe{(q_n,p_n)}$ of approximating solutions
obtained by solving just linear systems of ordinary differential equations.
Namely, by performing the above estimates 
on $(q_n,p_n)$ exactly in the same way as we did,
and using standard compactness results,
one finds a weak limit $(q,p)$ in the topologies associated to \eqref{formala}, \eqref{formala2}
and it is immediately clear that \pier{$(q,p)$} is a variational solution 
of the problem we want to solve.
Hence, Theorem~\ref{Existenceadj} actually holds.}\qed

At this point, we are ready to prove Theorem~\ref{CNoptadj} on optimality,
i.e., the necessary condition \eqref{cnoptadj} for
$\uopt$ to be an optimal control in terms of the solution \pier{$(q,p)$}
of the adjoint problem~\Pbladj.
So, we fix an arbitrary $u\in\Uad$
and use  the variational formulations of both the \lineariz ed problem 
(corresponding to $h=u-\uopt$) and the adjoint problem.

We test \eqref{linprima} by  $q$, \eqref{linseconda} by $p$, use \eqref{linterza}, and we take $z=-\Theta$ in \eqref{primaadj} and $w=-\Phi$ in \eqref{secondaadj}, respectively. \pier{Then, we
add all the equalities} we obtain to each other.
Most of the terms cancel out and \pier{we infer that}
\begin{align}
 & \intQ \kappa_1\Theta(\thetaopt-\theta_Q)+\intQ\kappa_2\Phi(\phiopt-\phi_Q)+\int_\Omega\kappa_3\Theta(T)(\thetaopt(T)-\theta_\Omega)+\int_\Omega\kappa_4(\phiopt(T)-\phi_\Omega)\non\\
 &=\intQ (u-u^*)q\geq 0\,.
  \non
\end{align}
As $u\in\Uad$ is arbitrary, this implies the pointwise inequality
 \eqref{cnoptadj} and the proof of Theorem~\ref{CNoptadj} is complete.\qed

%%%%%%%%%%%%%%%%%%%%%%%%%%%%%%%%%%%%%%%%%%%%%%%%%%%%%%%%%%%%%%%%%%%%%%%%

\section*{Acknowledgements}
This research activity has been performed in the framework of an
Italian-Romanian  three-year project on ``Control and 
stabilization problems for phase field and biological systems'' financed by the Italian CNR and the Romanian Academy.
Moreover, the financial support of  the project Fondazione Cariplo-Regione Lombardia  MEGAsTAR 
``Matema\-tica d'Eccellenza in biologia ed ingegneria come accelleratore 
di una nuova strateGia per l'ATtRattivit\`a dell'ateneo pavese'' is gratefully acknowledged by the authors. 
The present paper 
also benefits from the support of the MIUR-PRIN Grant 2015PA5MP7 ``Calculus of Variations'' for PC and GG, 
the GNAMPA (Gruppo Nazionale per l'Analisi Matematica, la Probabilit\`a e le loro Applicazioni) 
of INdAM (Istituto Nazionale di Alta Matematica) for PC, GG and~ER,
and \gabri{by a grant of Ministery of Research and Innovation, CNCS – UEFISCDI, 
project number PN-III-P4-ID-PCE-2016-0372, within PNCDI III} for GM.

%%%%%%%%%%%%%%%%%%%%%%%%%%%%%%%%%
%% bibliography
%%%%%%%%%%%%%%%%%%%%%%%%%%%%%%%%%

\vspace{3truemm}

%{\color{red}I JUST COPIED THE REFERENCE LIST  OF \ CGS \ AND ADDED WHAT I NEEDED.
%HENCE, THIS LIST HAS TO BE COMPLETELY REVISED !!!}

\Begin{thebibliography}{10}

\bibitem{Barbu}
V. Barbu,
``Nonlinear semigroups and differential equations in Banach spaces'',
\pier{Noordhoff,} 
Leyden, 
1976.

\bibitem{BBCG}
V. Barbu, M.L. Bernardi, P. Colli, G. Gilardi,
Optimal control problems of phase relaxation models,
J. Optim. Theory Appl. {\bf 109} (2001), 557--585.

\bibitem{BCF} 
J.L. Boldrini, B.M.C. Caretta, E. Fern{\'a}ndez-Cara,
Some optimal control problems for a two-phase field model of solidification,
Rev. Mat. Complut. {\bf 23} (2010), 49--75.

\bibitem{Brezis}
H. Brezis,
``Op\'erateurs maximaux monotones et semi-groupes de contractions
dans les espaces de Hilbert'',
North-Holland Math. Stud.
{\bf 5},
North-Holland,
Amsterdam,
1973.

\bibitem{BroHilNC}
D. Brochet, D. Hilhorst, A. Novick-Cohen, 
Maximal attractor and inertial sets for a
conserved phase field model, 
Adv. Differential Equations {\bf 1} (1996), 547--578.

\bibitem{BrokSpr} 
M. Brokate, J. Sprekels,
``Hysteresis and \pier{phase transitions}'',
Springer, New York, 1996.

\bibitem{Cag}
G. Caginalp,
The dynamics of a conserved phase field system: Stefan-like, Hele-Shaw, and Cahn-Hilliard models as asymptotic limits, 
IMA J. Appl. Math. {\bf 44} (1990), 77--94.

\bibitem{CH58}
J.W. Cahn and J.E. Hilliard, 
Free energy of a nonuniform system I. Interfacial free
energy, J. Chem. Phys. {\bf 2} (1958), 258--267.

\bibitem{CGW}
C. Cavaterra, M. Grasselli, H. Wu,
Non-isothermal viscous {C}ahn-{H}illiard equation with
inertial term and dynamic boundary \pier{conditions,
Comm. Pure Appl. Anal. {\bf 13} (2014), 1855--1890.} 

\bibitem{CFGS}
P. Colli, M.H. Farshbaf-Shaker, G. Gilardi, J. Sprekels, 
Optimal boundary control of a viscous Cahn\pier{--}Hilliard system with dynamic boundary condition and double obstacle potentials, 
SIAM J. Control Optim. {\bf 53} (2015), 2696--2721. 

\bibitem{CGLN}
P. Colli, G. Gilardi, P. Lauren\c cot, A. Novick-Cohen,
Uniqueness and long-time behaviour for the conserved phase-field  system with memory,
Discrete Contin. Dynam. Systems {\bf 5} (1999), 375--390.

\bibitem{CGM}
P. Colli, G. Gilardi, G. Marinoschi, 
A boundary control problem for a possibly singular phase field system
with dynamic boundary conditions,
J. Math. Anal. Appl. {\bf 434} (2016), 432--463.

\bibitem{CGMR1} 
P. Colli, G. Gilardi, G. Marinoschi, E. Rocca, 
Optimal control  for a phase field system with a possibly singular potential, 
Math. Control Relat. Fields {\bf 6} (2016), 95--112.

\pier{\bibitem{CGMR2} 
P. Colli, G. Gilardi, G. Marinoschi, E. Rocca, 
Distributed optimal control problems for 
phase field systems with singular potential, 
An. \c Stiin\c t. Univ. ``Ovidius'' Constan\c ta 
Ser. Mat., to appear (2017)}\betti{.}

\bibitem{CGS1}
P. Colli, G. Gilardi, J. Sprekels,
On the Cahn--Hilliard equation with dynamic boundary conditions
and a dominating boundary potential, J. Math. Anal. Appl. 
{\bf 419} (2014), 972--994.

\bibitem{CGS15}
P. Colli, G. Gilardi, J. Sprekels,
A boundary control problem for the pure Cahn-Hilliard equation with dynamic boundary conditions, 
Adv. Nonlinear Anal. {\bf 4} (2015), 311--325.

\bibitem{CGS16}
P. Colli, G. Gilardi, J. Sprekels, 
A boundary control problem for the viscous Cahn\pier{--}Hilliard equation with dynamic boundary conditions, 
Appl. Math. Optim. {\bf 73} (2016), 195--225. 

\bibitem{CGPS}
P. Colli, G.~Gilardi, P. Podio-Guidugli, J. Sprekels, 
Distributed optimal control of a nonstandard system of phase field equations,
Contin. Mech. Thermodyn. {\bf 24} (2012), 437--459. 

\bibitem{CGS} P. Colli, G. Gilardi, J. Sprekels, 
Analysis and optimal boundary control of a nonstandard system of phase field equations, 
Milan J. Math. {\bf 80} (2012), 119--149.

\bibitem{CMR}
P. Colli, G. Marinoschi, E. Rocca, 
Sharp interface control in a Penrose-Fife model, 
\pier{ESAIM Control Optim. Calc. Var.} {\bf 22} (2016), 473--499.

\bibitem{Gil}
 G. Gilardi, 
 On a conserved phase field model with irregular potentials and dynamic boundary conditions, 
 Rend. Cl. Sci. Mat. Nat. {\bf 141} (2007), 129--161. 

\bibitem{HoffJiang}
K.-H.Hoffmann, L.S. Jiang, 
Optimal control of a phase field model for solidification,
Numer. Funct. Anal. Optim. {\bf 13} (1992), 11--27. 

\bibitem{HKKY}
K.-H. Hoffmann, N. Kenmochi, M. Kubo, N. Yamazaki, 
Optimal control problems for models of phase-field type with hysteresis of play operator,
Adv. Math. Sci. Appl. {\bf 17} (2007), 305--336.

\bibitem{KenmNiez}
N. Kenmochi, M. Niezg\'odka, 
Nonlinear system for non–isothermal diffusive phase separation, 
J. Math. Anal. Appl. {\bf 188} (1994), 651--679.

\bibitem{LSU}
O.A. Lady\v zenskaja, V.A. Solonnikov, N.N. Ural'ceva:
``Linear and quasilinear equations of parabolic type'',
Trans. Amer. Math. Soc.~{\bf 23},
Amer. Math. Soc., Providence, RI,
1968.

\bibitem{LK} C. Lefter, J. Sprekels, 
Optimal boundary control of a phase field 
system modeling nonisothermal phase transitions,
Adv. Math. Sci. Appl. {\bf 17} (2007), 181--194.
		
\bibitem{Lions}
\pier{J.-L.~Lions,
``Quelques m\'ethodes de r\'esolution des probl\`emes
aux limites non lin\'eaires'',
Dunod; Gauthier-Villars, Paris, 1969.}

\bibitem{Mir}
A. Miranville, 
On the conserved phase-field model,
J. Math. Anal. Appl. {\bf 400} (2013), 143--152. 

\bibitem{MZ}
A. Miranville, S. Zelik, Robust exponential attractors for Cahn-Hilliard type equations
with singular potentials, Math. Methods Appl. Sci. {\bf 27} (2004)\pier{,} 545--582.

\bibitem{SY} 
K. Shirakawa, N. Yamazaki, Optimal control problems of phase field 
system with total variation functional as the interfacial energy,
Adv. Differential Equations {\bf 18} (2013), 309--350.

\bibitem{Simon}
J. Simon,
{Compact sets in the space $L^p(0,T; B)$},
{Ann. Mat. Pura Appl.~(4)\/} 
{\bf 146} (1987), 65--96.

\bibitem{SprZheng}
J. Sprekels, S. Zheng, 
Optimal control problems for a thermodynamically consistent model of phase-field type for phase transitions, 
Adv. Math. Sci. Appl. {\bf 1} (1992), 113--125.

\End{thebibliography}

\End{document}

\bye